\documentclass[final]{siamltex}

\usepackage{amssymb,amsmath}
\usepackage{amsfonts}
\usepackage[english]{babel}
\usepackage[latin1]{inputenc}
\usepackage{fancyhdr}
\usepackage{yfonts}
\usepackage{mathrsfs}
\usepackage[dvips]{graphicx}
\usepackage[rflt]{floatflt}
\usepackage{cite}

\usepackage{algorithm}
\floatname{algorithm}{Alg.}
\usepackage{algorithmic}
\usepackage{psfrag}

\providecommand{\com[2]}{\begin{tt}[#1: #2]\end{tt}}
\usepackage{color}

\newtheorem{prop}{Proposition}
\newtheorem{thm}{Theorem}

\newtheorem{cor}{Corollary}
\newtheorem{conj}{Conjecture}
\newtheorem{lem}{Lemma}

\providecommand{\prt}[1]{\left( #1 \right)}

\providecommand{\abs[1]}{\left|#1\right|}
\providecommand{\norm[1]}{\abs{\abs{#1}}}%

 \providecommand{\quotes[1]}{``#1"}

\usepackage{color}

\def\jnt#1{{#1}}
\def\lmod#1{{#1}}
\def\blue#1{{#1}}
\def\bo#1{{#1}}

\def\brn{\langle n \rangle}

\usepackage{epsfig}

\title{\jnt{Continuous-time average-preserving opinion dynamics with opinion-dependent communications}}

\author{Vincent D.
Blondel, Julien M. Hendrickx and John N. Tsitsiklis
\thanks{
This research was supported by the National Science Foundation
under grant ECCS-0701623, by the Concerted Research Action (ARC)
\quotes{Large Graphs and Networks} of the French Community of
Belgium, and by the Belgian Programme on Interuniversity
Attraction Poles initiated by the Belgian Federal Science Policy
Office. The scientific responsibility rests with its authors.
Julien Hendrickx holds postdoctoral fellowships from the
F.R.S.-FNRS (Belgian Fund for Scientific Research) and the
B.A.E.F. (Belgian American Education
Foundation). %
V.\ D.\ Blondel is with Department of Mathematical Engineering,
Universit\'e catholique de Louvain, Avenue Georges Lemaitre 4,
B-1348 Louvain-la-Neuve, Belgium; {\tt\small
vincent.blondel@uclouvain.be} J.\ M.\ Hendrickx  and J.\ N.\
Tsitsiklis are with the Laboratory for Information and Decision
Systems, Massachusetts Institute of Technology, Cambridge, MA
02139, USA; {\tt\small jm\_hend@mit.edu, jnt@mit.edu}. Part of
this reasearch was conducted when J. M. Hendrickx was at the
Université catholique de Louvain.}}

\begin{document}
\maketitle \thispagestyle{empty}
 \begin{abstract}
We study a simple continuous-time multi-agent \jnt{system} related to
Krause's model of opinion dynamics: \blue{each} agent holds a real value,
and this value is continuously attracted by \jnt{every other value}
differing from it by less than 1, with an intensity proportional
to the \jnt{difference.}

We prove convergence to a set of clusters, \jnt{with the agents in
each cluster sharing  a common} value, and provide a lower bound
on the distance between clusters \jnt{at a stable} equilibrium,
\jnt{under a suitable notion of multi-agent system stability.}

To better understand the behavior of the system for \jnt{a} large
number of agents, we introduce a variant involving a continuum of
agents. We prove, under some conditions, the existence of a
solution to \jnt{the system dynamics}, convergence to clusters,
and a non-trivial lower bound on the distance between clusters.
Finally, we establish that the \jnt{continuum model accurately}
represents the asymptotic behavior of \jnt{a system with a finite
but}  large number of agents.

\end{abstract}

\section{Introduction}

\begin{figure}
\centering
\includegraphics[scale=.5]{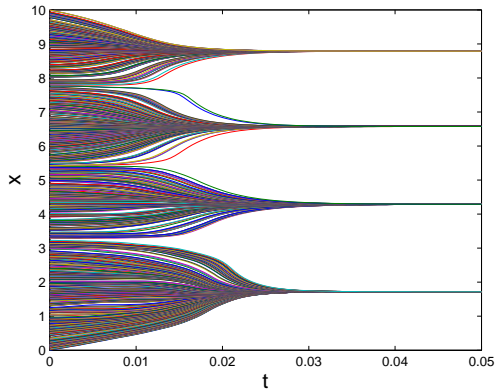}
\caption{Evolution with time of the values $x_i(t)$ for 1000
agents, with initial values \jnt{randomly and uniformly} distributed on
$[0,10]$. Observe the convergence to 4 clusters separated by
slightly more than 2.}\label{fig:ex_APCT}
\end{figure}

We study a continuous-time multi-agent model: \jnt{each of} $n$
agents, labeled $1,\dots,n$, maintains a real number
\jnt{(``opinion'')} $x_i(t)$, \lmod{which is a continuous function
of time and evolves according to the integral \bo{equation}
version of}
\begin{equation}\label{eq:def_system_intro_diff}
\dot x_i(t)= \sum_{j:\,\abs{x_i(t)-x_j(t)}<1}\prt{x_j(t) -x_i(t)}.
\end{equation}
This model has an interpretation in terms of opinion dynamics: an
agent considers \jnt{another} agent \jnt{to be a neighbor} if
their opinions differ by less than 1, and agent opinions are
continuously attracted by  \jnt{their neighbors' opinions.}
Numerical simulations show that the system converges to clusters
\jnt{inside}  which all agents share \jnt{a common} value.
\jnt{Different clusters lie at a distance of at least 1 from each
other}, and often approximately 2, as shown in Figure
\ref{fig:ex_APCT}. We focus on understanding these convergence
properties and the structure of the set of clusters, including
\jnt{the} asymptotic behavior for large $n$.

Observe that \jnt{the agent interaction topology in}
(\ref{eq:def_system_intro}) \emph{explicitly depends on the agent
states}, as $x_j(t)$ influences $x_i(t+1)$ only if
$\abs{x_i(t)-x_j(t)}<1$. Many multi-agent systems involve a
\jnt{changing} interaction topology; see e.g.
\cite{Tsitsiklis:84phdthesis,JadbabaieLinMorse:2003,Moreau:2005,
BlondelHendrickxOlshevskyTsitsiklis:2005,HendrickxBlondel:2006_MTNS,
VicsekCzirolBenjacobCohenSchchet:1995}, and
\cite{OlfatiSaberFaxMurray:2006, RenBeardAtkins:2007} for surveys.
In some cases, the interaction topology evolves randomly or
according to some exogenous scheme, but \jnt{in other cases it is
modeled as} a function of the agent states. With some exceptions
\cite{CuckerSmale:2005, CuckerSmale:2007,
JusthKrishnaprasad:2004}, \jnt{however,} this
\jnt{state-}dependence is not taken into account \jnt{in the
analysis,} probably due to the \jnt{technical difficulties that it
presents.}

\jnt{To address this issue,} we have \jnt{recently} analyzed
\cite{BlondelHendrickxTsitsiklis:2009_Krausemodel} one of the
simplest discrete-time multi-agent systems \jnt{with
state-dependent interaction topologies,} \jnt{namely,} Krause's
model\footnote{The model is sometimes referred to as \jnt{the}
Hegselmann-Krause model.} of opinion dynamics \cite{Krause:1997}:
$n$ agents maintain real numbers \jnt{(\quotes{opinions})}
$x_i(t)$, $i=1,\dots,n$, and synchronously update them \jnt{as
follows:}
\begin{equation*}
x_i(t+1) =
\frac{\sum_{j:\,\abs{x_i(t)-x_j(t)}<1}x_j(t)}{\sum_{j:\,\abs{x_i(t)-x_j(t)}<1}1}.
\end{equation*}
This model was particularly appealing due to its simple
formulation, and due to some peculiar \jnt{behaviors that it
exhibits}, which cannot be explained without taking into account
the explicit dynamics of the interaction topology. Indeed, a first
analysis using results on infinite inhomogeneous matrix products,
as in \cite{Lorenz:2005, HendrickxBlondel:2006_MTNS}, shows
\jnt{convergence} to clusters in which all agents share the same
opinion, and that the distance between any two clusters is at
least 1. Numerical simulations, however, show \jnt{a qualitative
behavior similar to the one shown in Figure \ref{fig:ex_APCT} for
the model (\ref{eq:def_system_intro}):} the distance between
\jnt{consecutive} clusters is usually significantly larger than 1,
and typically close to 2 when the number of agents is sufficiently
large, a phenomenon for which no explanation was available.

Our goal in \cite{BlondelHendrickxTsitsiklis:2009_Krausemodel} was
thus to develop a deeper understanding of Krause's model and of
these observed phenomena, by using explicitly the dynamics of the
interaction topology. \jnt{To this effect, we} introduced a new
notion of stability, \jnt{tailored to such} multi-agent systems,
\jnt{which provided} an explanation for the observed inter-cluster
distances when the number of agents is large. \jnt{Furthermore,}
to understand the asymptotic behavior \jnt{as the number of agents
increases, we also studied a model involving a} continuum of
agents. We obtained partial convergence results for this
\jnt{continuum model}, and proved nontrivial lower bounds on the
inter-cluster distances, under some conditions.

Our results in \cite{BlondelHendrickxTsitsiklis:2009_Krausemodel}
\jnt{were however incomplete in certain respects.}
In particular, the \jnt{question of} convergence
of the \jnt{continuum model} remains open, and \jnt{some of the} results
involve assumptions that are not \jnt{easy to check a priori}. We
see two main reasons for these difficulties. First, the system is
asymmetric, in the sense that the influence of $x_j(t)$ on
$x_i(t+1)$ can be very different from that on $x_i(t)$ on
$x_j(t+1)$, \jnt{when} $i$ and $j$ do not have the same number of
neighbors. Second, the discrete time nature of the system allows,
\jnt{for the continuum model,
buildup of an infinite concentration of agents with the same opinion,
thus}
breaking the continuity of the agent distribution.

\jnt{For the above reasons, we have} chosen to analyze here the
system (\ref{eq:def_system_intro}), a con\-ti\-nuous-time
symmetric \jnt{variant} of Krause's model, for which we provide
\jnt{crisper and more complete} results. \jnt{One reason is} that,
thanks to the symmetry, the average value $\frac{1}{n}\sum_i
x_i(t)$ is preserved, and the average value of a group of agents
evolves independent of the interactions taking place within the
group, unlike Krause's model. \jnt{In addition}, when two agent
values approach each other, their relative \jnt{velocity decays}
to zero, preventing the formation of infinite concentration in
finite time. The continuous-time nature of the system brings
\jnt{up} however some new mathematical challenges, related for
example to the existence and uniqueness of \jnt{solutions.}

\subsection{Outline and contributions}
In Section \ref{sec:discr_agents}, we give some basic properties
of \jnt{the model} (\ref{eq:def_system_intro}), and prove
convergence to clusters in which all agents share the same value.
We then analyze the distance between \jnt{consecutive clusters
building on an appropriate notion of stability with respect to
perturbing agents, introduced in}
\cite{BlondelHendrickxTsitsiklis:2009_Krausemodel}. This analysis
leads to a necessary and sufficient condition for stability that
is consistent with the experimentally observed inter-cluster
distances, and \jnt{to} a conjecture that the probability of
\jnt{convergence} to a stable equilibrium tends to one \jnt{as}
the number of agents \jnt{increases}.
In Section \ref{sec:continuum}, we introduce a variant involving a
continuum of agents, to approximate the \jnt{model for the case of
a finite but large number of agents.} Under some smoothness
assumptions on the initial conditions, we prove the existence of a
unique solution, convergence to clusters, and nontrivial lower
bounds on the inter-cluster distances, \jnt{consistent with the
necessary and sufficient for stability in the discrete-agent
model.}
Finally,  in Section \ref{sec:link_discrete-continuous}, we explore
the relation between \jnt{the} two models, and establish that the
behavior of the discrete model \jnt{approaches that of the
continuum model over finite but arbitrarily long time intervals,}
provided that the number of agents is sufficiently large.

The results summarized above differ from those those obtained
\jnt{in \cite{BlondelHendrickxTsitsiklis:2009_Krausemodel} for}
Krause's model, in three \jnt{respects}: \jnt{(i) we prove the
convergence of the continuum model, in contrast to the partial
results obtained for Krause's model; (ii) all of our stability and
approximation results are valid under some simple and easily
checkable smoothness assumptions on the initial conditions, unlike
the corresponding results in
\cite{BlondelHendrickxTsitsiklis:2009_Krausemodel} which require,
for example, the distance between the largest and smallest
opinions to remain larger than 2 at all times; (iii) finally, we
settle the problem of existence and uniqueness of a solution to
our equations, a problem that did not arise for Krause's
discrete-time model.}

\subsection{Related work}

Our model (\ref{eq:def_system_intro}) is \jnt{closely related} to that treated by
Canuto et al. \cite{CanutoFagnaniTilli:2008} who consider
multi-dimensional opinions whose evolution \jnt{is} described by
$$\dot x_i(t) = \sum_{j}
\xi\prt{x_i(t)-x_j(t)}\prt{x_i(t)-x_j(t)},$$ or, to a first order
discrete time approximation, \jnt{by}  $ x_i(t+\delta t) = x_i(t)
+ \delta t \dot x_i(t)$, where $\xi$ is a continuous\footnote{The
continuity assumption appears however unnecessary in the
discrete-time case, as was recently confirmed by \jnt{one of the
authors of \cite{CanutoFagnaniTilli:2008}} in a personal
communication.} nonnegative radially symmetric and decaying
function, taking positive values only for entries with norm
smaller than a certain \jnt{constant} $R$. Our model is
\jnt{therefore} a particular case of their continuous-time model
in one dimension with a step function for $\xi$, \jnt{except that}
a step function does does not satisfy their continuity assumption.

\jnt{The authors of \cite{CanutoFagnaniTilli:2008}} prove
\jnt{convergence of the opinions, in  distribution,} to clusters
separated by at least $R$ for \jnt{both discrete and
continuum-\blue{time} models.} Their convergence proof relies on
the decrease of the measured variance of the opinion distribution,
and is \jnt{based on} an Eulerian representation \jnt{that}
follows the density of agent \jnt{opinions}, \jnt{in contrast to
the Lagrangian representation  used in this paper, which} follows
the \jnt{opinion} $x$ of each agent. It is interesting to note
that despite the difference between \jnt{these two methods for
proving} convergence, \jnt{they} both appear to fail in the
absence of symmetry, and \jnt{cannot be used to prove} convergence
\jnt{for the continuum-agent variant of} Krause's model.

Finally, \jnt{the models in this paper} are also related to other classes of
rendezvous methods and opinion dynamics models, as described in
\cite{BlondelHendrickxTsitsiklis:2009_Krausemodel, Lorenz:2007}
and the references therein.

\section{Discrete agents}\label{sec:discr_agents}

\lmod{The differential \bo{equation}
(\ref{eq:def_system_intro_diff}) \bo{usually has no differentiable
solutions. Indeed,} observe that the right-hand side of the
equation can be discontinuous when the interaction topology
changes, which can prevent $x$ from being differentiable. To avoid
this difficulty, we consider functions $x:\Re^+\to \Re^n$ that are
solutions of the integral version of
(\ref{eq:def_system_intro_diff}),} namely
\begin{equation}\label{eq:def_system_intro}
x_i(t)= x_i(0) + \int_0^t
\sum_{j:\,\abs{x_i(\tau)-x_j(\tau)}<1}\prt{x_j(\tau)
-x_i(\tau)}d\tau.
\end{equation}
\lmod{Observe however that for all $t$ at which $\dot
x_i(t)$ exists, it can be computed using (\ref{eq:def_system_intro_diff}).}

\subsection{Existence and convergence}\label{sec:conv_and_basic_properties}

\jnt{Time-switched} linear systems are of the form
\lmod{$\bo{x}(t) = x(0) + \int_{0}^tA_\tau x(\tau)\, d\tau$,}
where $A_t$ is a piecewise constant function \jnt{of} $t$. They
always admit a unique solution provided that the number of
switches taking place during any finite time interval is finite.
\jnt{Position-switched} systems of the form \lmod{${\dot x}(t) =
x(0) + \int_{0}^tA_{x(\tau)} x(\tau)\, d\tau$} may on the other
hand admit none or multiple \jnt{solutions.} Our model
(\ref{eq:def_system_intro}) belongs to the latter class, and
\jnt{indeed admits} multiple solutions for some initial
conditions. Observe for example that the two-agent system with
initial condition $\tilde x = (-\frac{1}{2},\frac{1}{2})$ admits a
first solution $x(t) = \tilde x$, and a second solution $x(t) =
\tilde x e^{-t}$. The latter solution satisfies indeed the
\lmod{differential \bo{equation} (\ref{eq:def_system_intro_diff})
at every time except 0, and thus satisfies
(\ref{eq:def_system_intro})}. We will see however that such cases
are exceptional.

We say that $\tilde x \in \Re^n$ is a \emph{proper initial
condition of (\ref{eq:def_system_intro})} if:\\
(a) There exists a unique $x:\Re^+\to \Re^n:t\to x(t)$ satisfying
(\ref{eq:def_system_intro}), and such that $x(0)
= \tilde x$.\\
(b) The subset of $\Re^+$ on which $x$ is not differentiable is at
most countable, and \jnt{has no}  accumulation points.\\%
(c) If $x_i(t)=x_j(t)$ holds for some $t$, \bo{then $x_i(t')=x_j(t')$, for every $t'\geq t$.}

We then say that the solution $x$ is a \emph{proper solution} of
(\ref{eq:def_system_intro}). The \jnt{proof of the following
result is sketched} in Appendix~\ref{appen:exist_sol_discrete},
\lmod{and a detailed version is available in
\cite{BlondelHendrickxTsitsiklis:2009_proof_ex_unique_discrete}.}
\begin{thm}\label{thm:exist_unique_discr_agents}
Almost all $\tilde x\in \Re^n$ \jnt{(in the sense of Lebesgue measure)} are proper initial conditions.
\end{thm}

It follows from condition (c) and from the continuity of
proper solutions that if $x_i(t)\geq x_j(t)$ holds for some $t$,
then this inequality holds for all \jnt{subsequent times.} For the sake of clarity, we
assume thus in the sequel that the components of proper
initial conditions are sorted, that is,
\jnt{if $i
> j$, then} $\tilde x_i\geq \tilde
x_j$, which \jnt{also} implies that $x_i(t) \geq x_j(t)$ \jnt{for
all $t$.} Moreover, an explicit computation, which we perform in
Section \ref{sec:continuum} for a more complex system, shows that
$\abs{x_i(t)-x_j(t)} \geq \abs{\tilde x_i -\tilde x_j}e^{-nt}$.
Observe finally that if $x_{i+1}(t^*) -x_i(t^*) > 1$ holds for
some $t^*$ for a proper solution $x$, then $\dot x_{i+1}(t)\geq 0$
and $\dot x_{i}(t)\leq 0$ hold for \jnt{almost} all subsequent
$t$, so that $x_{i+1}(t) -x_i(t)$ remains larger\footnote{The case
$x_{i+1}(t) -x_i(t) =1$ is more complex. Agents could indeed
become \quotes{reconnected}, \jnt{as in the nonuniqueness example
given above} because (\ref{eq:def_system_intro}) \jnt{is allowed
to fail} at countably many times.} than 1. As a consequence, the
system can then be decomposed into \jnt{two} independent
subsystems, consisting of agents $1,\dots,i$, and $i+1,\dots,n$,
\jnt{respectively.}

We now characterize the evolution of the average  and
variance (sum of squared differences from the average)
\jnt{of the opinions}.
For this purpose, we \jnt{let $F$ be the set} of vectors $\tilde s
\in \Re^n$ \jnt{such} that for all $i,j\in
\{1,\dots,n\}$, either $\tilde s_i = \tilde s_j$, or $\abs{\tilde
s_i-\tilde s_j}\geq 1$. \jnt{We refer to vectors in $F$ as {\em \blue{equilibria.}}}

\begin{prop}
Let $x$ be a proper solution of (\ref{eq:def_system_intro}). The
average opinion $\bar x(t) = \frac{1}{n}\sum_{i=1}^nx_i(t)$ is
constant. The sum of squared differences from the average,
$V(x(t)) = \sum_{i=1}^n \prt{x_i(t) - \bar x(t)}^2$, \jnt{is
nonincreasing. Furthermore,} \blue{with the exception of a
countable set of times, if $x_t\not \in F$ (respectively, $x_t\in
F$), then the derivative $(d V/dt)(x(t))$ is negative
(respectively, zero).}
\end{prop}
\begin{proof}
For all $t$, except possibly \jnt{for} countably many,
\begin{equation}\label{eq:deriv_average}
\frac{d}{dt}\bar x(t) = \frac{1}{n}\sum_i\dot x_i(t) =
\frac{1}{n}\sum_{(i,j):\, |x_i(t)-x_j(t)|<1} \prt{x_j(t)-x_i(t)} = 0.
\end{equation}
\jnt{Since $x(t)$ is continuous, this
implies that $\bar x(t)$ is constant.}

Observe now that, for all $t$ at which $x$ is differentiable
$\frac{d}{dt}V(x(t))$ equals
\begin{equation*}
\sum_{i=1}^n 2 \prt{x_i(t) - \bar x(t)} \dot x_i(t)= 2
\sum_{i=1}^n x_i(t) \dot x_i(t) = 2\sum_{i=1}^n \sum_{j:\,
|x_i(t)-x_j(t)|<1}x_i(t) \prt{ x_j(t)-x_i(t)},
\end{equation*}
where we have used the relation (\ref{eq:deriv_average}) twice,
and the definition  (\ref{eq:def_system_intro}). The
\jnt{right-hand side of this equality} can be rewritten as
\begin{equation*}
\sum_{i,j:\, |x_i(t)-x_j(t)|<1}x_i(t) \prt{ x_j(t)-x_i(t)}
+\sum_{j,i:\, |x_j(t)-x_i(t)|<1}x_j(t) \prt{ x_i(t)-x_j(t)},
\end{equation*}
so that
\begin{equation*}
\frac{d}{dt}V(x(t)) = - \sum_{i,j:\, |x_i(t)-x_j(t)|<1}
\prt{x_j(t)-x_i(t)} ^2.
\end{equation*}
The latter expression is negative if $x(t)\not \in F$ and
\jnt{zero} otherwise.
\end{proof}

\jnt{There are several} convergence proofs for the system
(\ref{eq:def_system_intro}). We present here a simple \jnt{one}, which
\jnt{highlights} the importance of the average preservation
and symmetry properties, \jnt{and extends nicely to the continuum model.} A proof relying on other properties and
that can be used in the absence of symmetry can be found in
\cite{Hendrickx:2008phdthesis}.

\begin{thm}\label{thm:conv_cont-time_discr-agent}
Every proper solution $x$ of (\ref{eq:def_system_intro}) converges
to a limit $x^*\in F$; \jnt{that is,} for any $i, j$, \jnt{if
$x_i^* \not= x_j^*$, then $|x_i^*-x_j^*|\geq 1$.}
\end{thm}
\begin{proof}
Observe that by symmetry, \jnt{the equality}
\begin{equation*}
\sum_{i=1}^k\,\sum_{j\leq
k,\, \abs{x_i(t)-x_j(t)}<1}\prt{x_j(t)-x_i(t)} = 0
\end{equation*}
holds for any $k$ and any $t$. Therefore, it follows from
\eqref{eq:def_system_intro} that \lmod{for all $t$ but possibly
countably many,}
\begin{equation}\label{eq:d/dt_x1+x2+..+xk}
\frac{d}{dt}\sum_{i=1}^kx_i(t) =  \sum_{i=1}^k\, \sum_{j>
k,\, \abs{x_i(t)-x_j(t)}<1}\prt{x_j(t)-x_i(t)},
\end{equation}
which is nonnegative because $j>k>i$ implies $x_j(t)-x_i(t) \geq
0$. Since $\sum_{i=1}^kx_i(t)$ is bounded, this implies that it
converges \jnt{monotonically}, for any $k$. It then follows that every
$x_i(t)$ converges to a limit \jnt{$x_i^*$.} We assume that $x_k^*\neq x_{k+1}^*$ and suppose, to
obtain a contradiction, that $x^*_{k+1}-\jnt{x_k^*} <1$.
Then, since every \jnt{term} $x_j(t)-x_i(t)$ \jnt{on the right-hand side of \eqref{eq:d/dt_x1+x2+..+xk} is}
nonnegative, \jnt{the derivative on the left-hand side is asymptotically} positive and bounded away
from 0, preventing the convergence of $\sum_{i=1}^kx_i(t)$.
\jnt{Therefore, $x^*_{k+1}-x_k^*\geq 1$.}
\end{proof}

\subsection{\jnt{Stable equilibria} and inter-cluster
distances}\label{sec:discr_stab}

By the term \emph{clusters}, we will mean  the limiting values to
which the agent opinions converge. With some abuse of terminology,
we also refer to a set of agents whose opinions converge to the
same value as a cluster. Theorem
\ref{thm:conv_cont-time_discr-agent} implies that clusters are
separated by at least 1.
\jnt{On the other hand, extensive numerical experiments indicate that
the distance between
adjacent clusters is typically significantly larger than one, and if
the clusters contain the same number of agents, usually close to 2.}
\begin{figure}
\centering
\includegraphics[scale= .6]{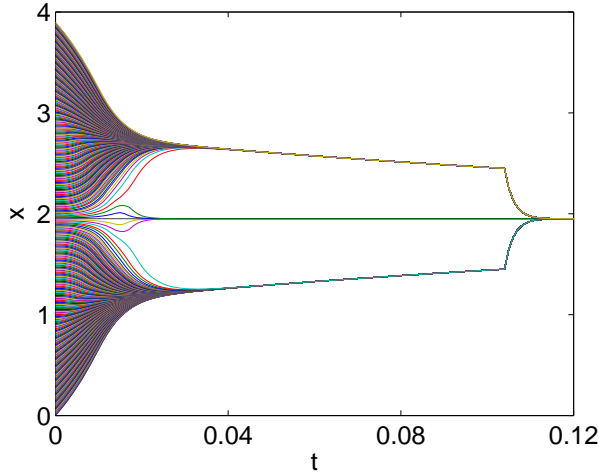}
\caption[Example of temporary \quotes{meta-stable}
equilibrium.]{Example of a temporary, \quotes{meta-stable,}
equilibrium. Initially, two clusters are formed and do not
interact with each other, but they both interact with a small
number of agents \jnt{in between.} As a result, the distance
separating them eventually becomes smaller than 1. The clusters
then attract each other directly and merge into \jnt{a single, larger}
cluster.} \label{fig:metastable_cont_time}
\end{figure}
We believe that this phenomenon can, at least partially, be
explained by the fact that clusters that are too close \jnt{to}
each other can be forced to merge by the presence of a small
number of agents between them, as in Figure
\ref{fig:metastable_cont_time}. To formalize this idea we
introduce a generalization of the system
(\ref{eq:def_system_intro}) in which each agent $i$ has a weight
$w_i$, and \jnt{its opinion} evolves according to
\begin{equation}\label{eq:def_weighted_system}
\lmod{x_i(t)= x_i(0) + \int_0^t
\sum_{j:\,\abs{x_i(\tau)-x_j(\tau)}<1}w_j\prt{x_j(\tau)
-x_i(\tau)}d\tau}.
\end{equation}
The results of Section \ref{sec:conv_and_basic_properties}
\jnt{carry over to the} weighted case \jnt{(the proof is the
same).} \jnt{We will refer to the sum of the weights of all agents
in a cluster, as its \emph{weight}.}  If all the agents \jnt{in} a
cluster have exactly the same \jnt{opinion,} the cluster behaves
as \jnt{a single} agent with this particular
weight\footnote{\lmod{In the case of non-proper initial conditions
leading to multiple solutions, there exists at least one solution
in which each cluster behaves as a single agent with the
corresponding weight.}}.

Let \jnt{$\tilde s\in F$ be an equilibrium vector.} Suppose that
\jnt{we} add a new agent of weight $\delta$ and \jnt{initial}
opinion \blue{$x_0$, consider the resulting configuration as an
initial condition, and let the system evolve according to some
solution $x(t)$ (we do not require uniqueness). We define
$\Delta(\delta,\tilde s)$ as the supremum of $\abs{x_i(t)-\tilde
s_i }$, where the supremum is taken over all possible initial
opinions $x_0$ of the perturbing agent, all $i$, all times $t$,
and all possible solutions $x(t)$ of the system
(\ref{eq:def_system_intro}). We say that $\tilde s$ is
\emph{stable} if $\lim_{\delta\downarrow 0} \Delta(\delta, \tilde
s)=0$.} An equilibrium is thus unstable if some modification of
fixed size can be achieved by adding an agent of arbitrarily small
weight. This notion of stability is \blue{almost} the same as the
one that we introduced for Krause's model in
\cite{BlondelHendrickxTsitsiklis:2007ECC,
BlondelHendrickxTsitsiklis:2009_Krausemodel}.

\begin{thm}\label{thm:stab_disc_agent}
An equilibrium is stable  if and only if for any two clusters $A$
and $B$ with weights $W_A$ and $W_B$, respectively, their distance
is  \blue{greater than}
$\jnt{d=}1+\frac{\min\left\{W_A,W_B\right\}}{\max\left\{W_A,W_B\right\}}$.
\end{thm}

\begin{proof}
The proof is very similar to the proof of Theorem 2 in
\cite{BlondelHendrickxTsitsiklis:2009_Krausemodel}. The main idea
is the following. A perturbing agent can initially be connected to
at most two clusters, and \jnt{cannot perturb the equilibrium
substantially} if it is connected to none or one. If it is
connected to two clusters $A,B$, it moves in the direction of
their center of mass $\frac{W_A\blue{\tilde s}_A+W_B\blue{\tilde
s}_B}{W_A+W_B}$, while the two clusters move at a much slower
pace, proportional \jnt{to} the perturbing agent's weight. \jnt{We
note that, by a simple algebraic calculation, the center of mass
of two clusters is within unit distance from both clusters if and
only if their distance is no more than $d$.}

\jnt{If the distance between the two clusters is more than $d$,
then} the center of mass of the two clusters is  more than
\jnt{unit distance away from one of the clusters, say from} $B$.
\jnt{Therefore, eventually} the perturbing agent is \jnt{no longer
connected to} $B$, and rapidly joins cluster $A$, having modified
the cluster positions only proportionally to its weight.
\jnt{Thus, the equilibrium is stable.}

On the other hand, \jnt{if the distance between the two clusters
is less than $d$, then} the center of mass is \jnt{less than unit
distance away}  from \jnt{both} clusters. \jnt{We can place the
perturbing agent at the center of mass. Then, the perturbing agent
does not move, but keeps attracting the two clusters, until
eventually they become connected and then rapidly merge. Thus, the
equilibrium is not stable.}

\blue{If the distance between clusters is exactly equal to $d$,
the center of mass is at exactly unit distance from one of the two
clusters. Placing a perturbing agent at the center of mass results
in nonunique solutions. In one of these solutions, the clusters
start moving towards their center of mass,} \lmod{and the
subsequent behavior is the same as in the case where the distance
between clusters is smaller than $d$, thus again showing
instability. Such a solution violates the differential
version of (\ref{eq:def_system_intro}) only at time $t=0$ and
thus satisfies (\ref{eq:def_system_intro}).}
\end{proof}

Theorem \ref{thm:stab_disc_agent} characterizes stable equilibria
in terms of a lower bound on inter-cluster distances. It allows
for inter-cluster distances at a stable equilibrium that are
smaller than 2, provided that the clusters have different weights.
This is consistent with experimental observations for certain
initial opinion distributions (see \cite{Hendrickx:2008phdthesis}
for example). On the other hand, for the frequently observed case
of clusters with equal weights, stability requires inter-cluster
distances \jnt{of} at least 2. Thus, this result comes close to a
full explanation of the observed inter-cluster distances of about
2.2. \jnt{Of course,} there is no guarantee that our system will
converge to a stable equilibrium. (A trivial example is obtained
by initializing the system at an unstable equilibrium.)
\lmod{However,} we have observed that for a given distribution of
initial opinions, and as the number of agents increases, we almost
always obtain convergence to a stable equilibrium. This leads us
to the following conjecture.

\begin{conj}\label{conj:stab_discr_agent}
Suppose that the initial opinions are chosen randomly and
independently according to a \blue{bounded probability density
function with connected support, which is also bounded below by a
positive number on its support.} Then, the probability of
convergence to a stable equilibrium tends to 1, as the number of
agents increases to infinity.
\end{conj}

In addition to extensive numerical \jnt{evidence} (see for example
\cite{Hendrickx:2008phdthesis}), this conjecture is supported by
the intuitive idea that if the number of agents is sufficiently
large, convergence to an unstable equilibrium is made impossible
by the presence of at least one agent connected to the two
clusters. It is also supported by results obtained in the next
sections. \jnt{A} similar conjecture has been made for Krause's
model \cite{BlondelHendrickxTsitsiklis:2007ECC,
BlondelHendrickxTsitsiklis:2009_Krausemodel}.

\section{Agent continuum}\label{sec:continuum}

To further analyze the properties of (\ref{eq:def_system_intro})
and its behavior as the number of agents increases, we now
consider a \jnt{variant involving} a
continuum of agents. We use the interval $I=[0,1]$ to index the
agents, and denote by \blue{$Y$} the set of bounded measurable functions
$\tilde x:I\to\Re$, attributing an opinion $\tilde x(\alpha) \in
\Re$ to every agent in $I$. As an example, a uniform distribution
of opinions is given by $\tilde x(\alpha) = \jnt{\alpha}$. We use the
function $x: I\times \Re^+ \to \Re :(\alpha,t)\to x_t(\alpha)$ to
\jnt{describe the collection of all opinions at different times.}\footnote{\jnt{Note the reversal of notational conventions: the subscript now indicates time rather than an agent's index.}} \jnt{We} denote by $x_t$ the function in \blue{$Y$} obtained
by restricting $x$ to a certain value of $t$. For a given initial
opinion function $\tilde x_0 \in \blue{Y}$, we are interested in
functions $x$ satisfying
\begin{equation}\label{eq:def_cont_agent_derivative}
\frac{d}{dt} x_t(\alpha) = \int_{\beta:\,(\alpha,\beta)\in C_{x_t}}
\prt{x_t(\beta)-x_t(\alpha)}d\beta,
\end{equation}
where $C_{\lmod{\tilde x}}\subseteq I^2$ is defined for any
$\lmod{\tilde x}\in \blue{Y}$ by
\begin{equation*}
C_{\lmod{\tilde x}}  := \{(\alpha,\beta)\in I^2:\abs{\lmod{\tilde
x}(\alpha)- \lmod{\tilde x}(\beta)}< 1\}.
\end{equation*}
In the sequel, we denote by $\chi_{\lmod{\tilde x} }$ the
indicator functions of $C_{\lmod{\tilde x}}$.

Note that $x_0$, the restriction of $x$ to $t=0$, should not be
confused with $\tilde x_0$, an arbitrary function in \blue{$Y$}
intended as an initial condition, but for which they may possibly
exist none or several corresponding functions $x$. The existence
or uniqueness of a solution to
\jnt{\eqref{eq:def_cont_agent_derivative}} is not guaranteed, and
there may moreover exist functions that satisfy this equation in a
weaker sense, \jnt{without being  differentiable in $t$. For this
reason, it is more convenient to formally define the model through
an integral equation.}  For an initial opinion function $\tilde
x_0\in \blue{Y}$, we \jnt{are interested in measurable functions}
$x: I\times \Re^+ \to \Re :(\alpha,t)\to x_t(\alpha)$ such that
\begin{equation}\label{eq:def_cont_agent_integral}
x_t(\alpha) = \tilde x_0(\alpha) + \int_{\jnt{0}}^t
\prt{\int_{\beta:\,(\alpha,\beta)\in C_{x_\tau}}
\prt{x_\tau(\beta)-x_\tau(\alpha)}\, d\beta}d\tau
\end{equation}
holds for every $t$ and for every $\alpha \in I$.\footnote{\jnt{A slightly more general definition would require
(\ref{eq:def_cont_agent_integral}) to be satisfied for almost all
$\alpha \in I$. However, this would result in distracting
technicalities.}}
\jnt{Similar to} the case of discrete agents, one can easily prove that for
any solution $x$ of (\ref{eq:def_cont_agent_integral}), $\bar x_t
:= \int_{0}^1 x_t(\alpha)\, d\alpha$ is constant, and $\int_{0}^1
\prt{x_t(\alpha)-\bar x_t}^2  d\alpha$ is nonincreasing in $t$.

For the sake of simplicity, we \jnt{will} restrict attention  to
nondecreasing (and often increasing) opinion functions, and
\jnt{define} $X$ as the set of nondecreasing bounded functions
$\tilde x :I\to \Re$. This is \jnt{no essential loss of
generality, because the only quantities of interest relate to the
distribution of opinions; furthermore, monotonicity of initial
opinion functions can be enforced using a measure-preserving
reindexing of the agents; finally, monotonicity is preserved by
the dynamics under mild conditions.} \jnt{In the sequel an element
of $X$ will be referred as a nondecreasing function. Furthermore,
if $x:I\times \blue{[0,\infty)}\to \Re$ is such that $x_t\in X$
for all $t$, we will also say that $x$ is nondecreasing.}

\subsection{Existence and uniqueness of \jnt{solutions}}

The existence of a unique solution to
(\ref{eq:def_cont_agent_integral}) is in general not guaranteed,
as there \jnt{exist} initial conditions allowing for multiple
solutions. Consider for example $\tilde x_0(\alpha) = -1/2$ if
$\alpha \in [0,\frac{1}{2}]$, and $\tilde x_0(\alpha) =1/2$
\jnt{otherwise.} Observe that, similar to our discrete-agent
example, \jnt{$x_t = \tilde x_0$ and $x_t(\alpha) = \tilde
x_0(\alpha)e^{-t}$ are two possible solutions of
(\ref{eq:def_cont_agent_integral}).} Nevertheless, we \jnt{will
prove existence and uniqueness} provided that the initial
condition, \jnt{as a function of $\alpha$,}  has a positive and
bounded increase rate; \jnt{this is equivalent to}  assuming that
the density of initial opinions is bounded from above and from
below on \jnt{its support, \lmod{which is connected.}}

Our proof of existence and uniqueness is based on
the Banach fixed point theorem, applied to the operator $G$
\jnt{that maps measurable functions $x:I\times
[0,t_1]\to \Re$ into the set of such functions, according to}
\begin{equation*}
(G(x))_t(\alpha) = \tilde x_0(\alpha) + \int_{\tau = 0}^t
\prt{\int_{\beta:\, (\alpha,\beta)\in C_{x_\tau}}
\prt{x_\tau(\beta)-x_\tau(\alpha)}\,d\beta}d\tau,
\end{equation*}
for some \jnt{fixed} $t_1$. Observe indeed that $x$ is a solution of the
system (\ref{eq:def_cont_agent_integral}) if and only if \jnt{$x_0=\tilde x$ and} $x=G(x)$.

It is convenient to introduce some \jnt{additional} notation. For
positive real numbers $m,M$, we call $X_m\subset X$ the set of
\jnt{nondecreasing}
functions $\tilde x: I\to\Re$ such that
\begin{equation*}
\frac{\tilde x(\beta)-\tilde x(\alpha)}{\beta-\alpha} \geq m
\end{equation*}
holds for every $\beta\not = \alpha$, and $X^M\subset X$ the set
of \jnt{nondecreasing} functions $\tilde x$ such that for all $\beta\not = \alpha$,
\begin{equation*}
\frac{\tilde x(\beta)-\tilde x(\alpha)}{\beta-\alpha} \leq M.
\end{equation*}
We then denote $X_m\cap X^M$ by $X_m^M$, and say that a function
$\tilde x \in X$ is \emph{regular} if it belongs to $X_m^M$ for
some $m,M>0$.
\providecommand{\LL}{\mathscr{L}}%
Let now $\LL$ be the operator defined on $X$ and taking its values
in the set of functions from $I$ to \jnt{$\Re$}, defined by
\begin{equation}\label{eq:def_operator_LL}
\LL(\tilde x)(\alpha) = \int \chi_{\tilde x}(\alpha,\gamma)
\prt{\tilde x(\gamma)-\tilde x(\alpha)} d\gamma.
\end{equation}
Observe that (\ref{eq:def_cont_agent_integral}) can be rewritten
as $x_t(\alpha) = \tilde x_0(\alpha) + \int_{0}^t
\LL(x_\tau)(\alpha)\, d\tau = (G(x))_t$.

\jnt{The proof of existence and uniqueness rests on two important
qualitative properties of our model. The first, given in Lemma
\ref{lem:lipscont} below, establishes that $\LL$ is Lipschitz
continuous on $X_m$. This property will allow us to establish that
the operator $G$ is a contraction (when $t_1$ is small enough),
and to apply Banach's fixed point theorem. The second, Lemma
\ref{lem:muderiv} below, gives bounds on the rate at which the
opinions of different agents can approach each other. It is
instrumental in showing that regularity is preserved, allowing us
to apply the same argument and extend the solution to arbitrarily
long time intervals.}

\begin{lem}\label{lem:lipscont}
Let $\tilde x$ be a function in $X_m$, where $m >0$. \jnt{The
operator} $\LL$ is \jnt{Lipschitz continuous \blue{at $\tilde x$}}
with respect to the $\norm{\,\cdot \,}_\infty$ norm. More
precisely, for any $\tilde y \blue{\in Y}$,
\begin{equation*}
\norm{\LL(\tilde x)-\LL(\tilde y)}_\infty \leq
\prt{2+\frac{8}{m}}\norm{\tilde x- \tilde y}_\infty.
\end{equation*}
\end{lem}

\begin{proof}
\jnt{Let \blue{$\tilde x\in X_m$, $\tilde y\in Y$,} and $\delta=
\norm{\tilde x-\tilde y}_\infty$. Fix some} $\alpha\in I$, and
\jnt{let} $N_x := \{\gamma: \abs{\tilde x (\gamma) -\tilde x
(\alpha) } < 1 \}$, $N_y := \{\gamma: \abs{\tilde y (\gamma)
-\tilde {\jnt{y}} (\alpha) } < 1 \}$ be the \jnt{sets} of agents
connected to $\alpha$, \jnt{under the configuration}
 defined by $\tilde x$ and $\tilde y$, respectively. Let \jnt{also}
$N_{xy} = N_x \cap N_y$, $N_{x\setminus y} = N_x \setminus
N_{xy}$, and $N_{y\setminus x} = N_y \setminus N_{xy}$. By the
definition (\ref{eq:def_operator_LL}) of $\LL$, we have
$\LL(\tilde x)(\alpha) = \int_{N_x}\prt{\tilde x (\gamma) - \tilde
x(\alpha)}d\gamma$ and $\LL(\tilde y)(\alpha) =
\int_{N_y}\prt{\tilde y (\gamma) - \tilde y(\alpha)}d\gamma$.
Therefore,
\begin{equation*}
\begin{array}{llll}
\LL(\tilde y)(\alpha) - \LL(\tilde x)(\alpha) &=&&
\int_{N_{xy}}\prt{\tilde y (\gamma) - \tilde x (\gamma) - \tilde
y(\alpha) + \tilde x(\alpha) } d\gamma \\ &&+&\int_{N_{y\setminus
x}} \prt{\tilde y (\gamma) -\tilde y(\alpha) } d\gamma-
\int_{N_{x\setminus y}} \prt{\tilde x (\gamma) -\tilde x(\alpha) }
d\gamma.
\end{array}
\end{equation*}
It follows from the definition of $N_x$ and $N_y$ that
$\abs{\tilde x (\gamma) -\tilde x(\alpha)}<1$ holds for every
$\gamma \in N_{x\setminus y}\subseteq N_x$ and $\abs{\tilde y
(\gamma) -\tilde y(\alpha)}<1$ holds for every $\gamma \in
N_{y\setminus x}\subseteq N_y$. This leads to
\begin{equation}\label{eq:bound_L(y)-L(x)}
\begin{array}{ll}
\abs{\LL(\tilde y)(\alpha) - \LL(\tilde x)(\alpha)} &\leq
\int_{N_{xy}}\prt{\abs{\tilde y (\gamma) - \tilde x (\gamma)}
+\abs{ \tilde y(\alpha) - \tilde x(\alpha)} }d\gamma +
|N_{x\setminus y}| + |N_{y\setminus x}|\\
&\leq 2 \abs{N_{xy}}\delta + \abs{N_{x\setminus y}} +
\abs{N_{y\setminus x}}\\ &\leq 2 \delta + \abs{N_{x\setminus y}} +
\abs{N_{y\setminus x}},
\end{array}
\end{equation}
where we have used the bound $\abs{N_{xy}} \jnt{\leq} \abs{I}=1$ to obtain
the last inequality. It remains to give bounds on
$\abs{N_{x\setminus y}}$ and $\abs{N_{y\setminus x}}$.

If $\jnt{\gamma \in } N_{y\setminus x}$, then \jnt{$\gamma\in
N_y$, and} $\abs{\tilde y(\gamma)-\tilde y(\alpha)}< 1$. This
implies that
\begin{equation*}\label{eq:bound_d+1+d}
\abs{\tilde x(\gamma)-\tilde x(\alpha)} \leq \abs{\tilde
x(\gamma)-\tilde y(\gamma)} + \abs{\tilde y(\gamma) -\tilde
y(\alpha)} + \abs{\tilde y(\alpha)-\tilde x(\alpha) }\leq \delta +
1 + \delta.
\end{equation*}
Since the same $\gamma$ does not belong to $N_x$, \jnt{we also have}
$\abs{\tilde x(\gamma)-\tilde x(\alpha)}\geq 1$.
\jnt{Thus, for every $\gamma\in N_{y\setminus x}$, the opinion $\tilde{x}(\gamma)$ lies in the set
$$\left[\tilde
x(\alpha)-1 \jnt{-2\delta} , \tilde x(\alpha)-1 \right] \cup
\left[\tilde x(\alpha)+1,\tilde x(\alpha)+1 \jnt{+2\delta}
\right],$$ which has length at most $4\delta$. Since the rate of
change of opinions (with respect to the index $\gamma$) is at
least $m$, we conclude that $\abs{N_{y\setminus x}}\leq
4\delta/m$.} A similar argument shows that $\abs{N_{x\setminus
y}}\leq 4\delta/m$. \jnt{The inequality (\ref{eq:bound_L(y)-L(x)})
then becomes}
\begin{equation*}
\abs{\LL(\tilde y)(\alpha) - \LL(\tilde x)(\alpha)} \leq 2\delta
 + 8\frac{\delta}{m} = \prt{2+\frac{8}{m}}\norm{\tilde y - \tilde x}_\infty,
\end{equation*}
which \jnt{is} the desired result. \end{proof}

\begin{lem}\label{lem:muderiv}
Let \jnt{$\tilde x \in \blue{Y}$. Suppose that $\alpha,\beta\in I$, and \blue{$x(\alpha)\leq x(\beta)$.} Then,}
\begin{equation*}
\LL(\tilde x)(\beta) - \LL(\tilde x)(\alpha) \geq - \prt{\tilde
x(\beta) - \tilde x(\alpha) }.
\end{equation*}
\jnt{Furthermore, if} $\tilde x \in X_m$ for some $m>0$, then \begin{equation*}
\LL(\tilde x)(\beta) - \LL(\tilde x)(\alpha) \leq \frac{2}{m}
\prt{\tilde x(\beta) - \tilde x(\alpha) }.
\end{equation*}
\end{lem}

\begin{proof}
Let $N_\alpha := \{\gamma: \abs{\tilde x (\gamma) -\tilde x
(\alpha) } < 1 \}$ and $N_\beta := \{\gamma: \abs{\tilde x
(\gamma) -\tilde x (\beta) } < 1 \}$ be the \jnt{sets} of agents
connected to $\alpha$ and $\beta$, respectively. Let now
$N_{\alpha\beta} = N_\alpha \cap N_\beta$, $N_{\alpha\setminus
\beta }= N_\alpha \setminus N_{\alpha\beta}$ and $N_{\beta
\setminus \alpha} = N_\beta \setminus N_{\alpha\beta}$. It follows
from the definition (\ref{eq:def_operator_LL}) of $\LL$
that
\begin{equation}\label{eq:lxb_lxa}
\begin{array}{lllllll}
\LL(\tilde x)(\beta) &= & \int_{N_{\alpha\beta}} \prt{\tilde
x(\gamma)-\tilde x(\beta)}d\gamma &+& \int_{N_{\beta\setminus
\alpha}} \prt{\tilde
x(\gamma)-\tilde x(\beta)}d\gamma,\\
\LL(\tilde x)(\alpha) &= &
 \int_{N_{\alpha\beta}} \prt{\tilde
x(\gamma)-\tilde x(\alpha)}d\gamma &+& \int_{N_{\alpha\setminus
\beta}} \prt{\tilde x(\gamma)-\tilde x(\alpha)}d\gamma.
\end{array}
\end{equation}
\blue{The definitions} of the sets $N_{\beta\setminus \alpha}$ and
$N_{\alpha\setminus \beta}$, \lmod{together with $\tilde
x(\beta)\geq \tilde x(\alpha)$,} imply that $\tilde x (\gamma)
> \tilde x(\alpha)$ holds for all $\gamma \in N_{\beta\setminus
\alpha}$, and  $\tilde x (\gamma) < \tilde x(\beta)$ holds for
every $\gamma \in N_{\alpha \setminus \beta}$. Using these
inequalities and subtracting the two \jnt{equalities} above, we
obtain
\begin{equation*}
\LL(\tilde x)(\beta) - \LL(\tilde x)(\alpha) \geq
\int_{N_{\alpha\beta}} \prt{\tilde x(\alpha)-\tilde
x(\beta)}d\gamma + \int_{N_{\beta\setminus \alpha}\cup
N_{\alpha\setminus \beta}} \prt{\tilde x(\alpha)-\tilde
x(\beta)}d\gamma
\end{equation*}
Since $\abs{N_{\alpha\beta}}+
\abs{N_{\alpha\setminus \beta}}+\abs{N_{\beta\setminus \alpha}} =
\abs{N_\alpha \cup N_\beta} \leq \abs{I}=1$, \jnt{we obtain the first part of the lemma.}

Let us now assume that $\tilde x \in X_m$. It follows from
(\ref{eq:lxb_lxa}) and from the inequality $\tilde x(\beta) \geq
\tilde x(\alpha)$ that
\begin{equation}\label{eq:lxb-lxa_upperbound}
\LL(\tilde x)(\beta) - \LL(\tilde x)(\alpha) \leq
\int_{N_{\beta\setminus \alpha}} \prt{\tilde x(\gamma)-\tilde
x(\beta)}d\gamma - \int_{N_{\alpha\setminus \beta}} \prt{\tilde
x(\gamma)-\tilde x(\alpha)}d\gamma,
\end{equation}
which is bounded by $\abs{N_{\beta\setminus \alpha}}+
\abs{N_{\alpha\setminus \beta}}$. Observe that $\tilde
x(N_{\beta\setminus \alpha}) \jnt{\subseteq} [\tilde x(\alpha) +1,
\tilde x(\beta) +1)$. \jnt{Since} $\tilde x \in X_m$, we have
\begin{equation*}
\abs{N_{\beta\setminus \alpha}} \leq \frac{1}{m} \abs{\tilde
x(N_{\beta\setminus \alpha})} = \frac{1}{m}\prt{\tilde x(\beta) -
\tilde x(\alpha)}.
\end{equation*}
The same bound holds on $\abs{N_{\alpha\setminus \beta}}$. The
second \jnt{part of the} lemma follows from the bound on
(\ref{eq:lxb-lxa_upperbound}).
\end{proof}

\jnt{The first part of Lemma \ref{lem:muderiv}} implies that for
\blue{any solution $x$,} the \jnt{difference of} the opinions of
two agents decreases at most exponentially \jnt{fast}. As a
consequence, if the initial \jnt{condition}  of $x$ is \jnt{an
increasing function of $\alpha$,} then $x_t$ is \jnt{also
increasing} for all $t$. \jnt{We note that this statement does not
necessarily hold for non-regular initial conditions initial
conditions that are only nondecreasing.}

We can now formally state our existence and uniqueness result,
\jnt{with the rest of the proof given} in
Appendix~\ref{appen:exit_unique_continuous}. \jnt{This result also
shows that if the initial condition is regular, then the two
models given by a differential or integral equation, respectively,
admit a unique and common solution, which is regular at all
times.}

\begin{thm}\label{thm:existence_continuous}
\jnt{Suppose that the initial opinion function satisfies $\tilde
x_0\in X_m^M$, for some $m,M >0$.} Then \jnt{the models}
(\ref{eq:def_cont_agent_derivative}) and
(\ref{eq:def_cont_agent_integral}) admit a unique and common
solution $x$, \jnt{and $x$} satisfies
\begin{equation}\label{eq:mt<inc_rate<Mt}
me^{-t}\leq \frac{x_t(\beta)-x_t(\alpha)}{\beta-\alpha} \leq
Me^{4t/m},
\end{equation}
for every $t$ and $\beta\not= \alpha$.
\end{thm}

\subsection{Convergence and fixed points}\label{sec:convergence_continuous}

In this section, we prove that opinions converge to clusters
separated by at least \jnt{unit distance}, as in the case of discrete agents. The
proof \jnt{has} some similarities with the one of Theorem
\ref{thm:conv_cont-time_discr-agent}. It \jnt{involves} three
partial results, the first of which establishes the convergence of
the average value of $x_t$ on any interval.
\jnt{Lemma \ref{l:tail} below involves an assumption that $x_t$ is nondecreasing. By Theorem
\ref{thm:existence_continuous}, this is guaranteed if the
initial condition is regular.}

\begin{lem}\label{l:tail}
Let $x$ be a \jnt{nondecreasing} solution of the integral equation
(\ref{eq:def_cont_agent_integral}).  For any $c\in I$, the limit
\begin{equation*}
\lim_{t\to \infty} \int_{\jnt{0}}^c x_t(\alpha)\, d\alpha
\end{equation*}
exists. As a result, the average value \jnt{$(\int_b^c x_t(\alpha)\, d\alpha)/(c-b)$} of $x_t$ on any \jnt{positive length} interval $[b,c]$
converges as $t\to \infty$.
\end{lem}
\begin{proof}
Fix some $c\in[0,1]$ and $t_1,t_2$ with $0\leq t_1<t_2$. The
evolution equation (\ref{eq:def_cont_agent_integral}) yields
\begin{equation}\label{eq:inta}
\int_{0}^c x_{t_2}(\alpha) d\alpha = \int_0^c x_{t_1}(\alpha)
d\alpha + \int_{t_1}^{t_2}\prt{\int_{0}^c\int_0^1
\chi_{x_\tau}(\alpha,\beta)\prt{x_\tau(\beta)-x_\tau(\alpha)} d\beta
\,d\alpha } d\tau,
\end{equation}
where we have used the Fubini theorem to interchange the
integration with respect to $\tau$ and $\alpha$. We observe that
$$\int_0^c\int_0^c \chi_{x_\tau}(\alpha,\beta) \prt{x_\tau(\beta)-x_\tau(\alpha) }
d\beta\, d\alpha=0,$$ because of the symmetry property
$\chi_{x_\tau}(\alpha,\beta) =  \chi_{x_\tau}(\beta,\alpha)$.
Therefore,

\begin{equation}\label{eq:intb}
\int_0^c\int_0^1 \chi_{x_\tau}(\alpha,\beta)
\prt{x_\tau(\beta)-x_\tau(\alpha)}\, d\beta\, d\alpha =
\int_0^c\int_c^1 \chi_{x_\tau}(\alpha,\beta)
\prt{x_\tau(\beta)-x_\tau(\alpha)}\, d\beta\, d\alpha
\end{equation}
The latter integral is nonnegative, because
$x_\tau(\beta)-x_\tau(\alpha)\geq 0$ whenever $\alpha\leq
c\leq\beta$. \jnt{Thus, $\int_0^c x_t(\alpha)\, d\alpha$ is a bounded and nondecreasing function of $t$, hence converges, which is} the desired result.
\end{proof}

\begin{prop}
Let $x$ be a solution of the integral equation
(\ref{eq:def_cont_agent_integral}) such that $x_t$ is
nondecreasing \jnt{in $\alpha$} for all $t$. For all $\alpha\in I$, except possibly
for a countable set, the limit $\lim_{t\to\infty}x_t(\alpha)$
exists.
\end{prop}
\begin{proof}
Let $\tilde y(\alpha)=\limsup_{t\to\infty} x_t(\alpha)$. Since,
for any $t$, $x_t(\alpha)$ is a nondecreasing function of
$\alpha$, it follows that $\tilde y(\alpha)$ is also nondecreasing
in $\alpha$. Let $S$ be the set of all $\alpha$ at which $\tilde
y(\cdot)$ is discontinuous. Since $\tilde y(\cdot)$ is
nondecreasing, it follows that $S$ is at most countable.

Fix some $\alpha\notin S$. Suppose, in order to derive a
contradiction, that $x_t(\alpha)$ does not converge to $\tilde
y(\alpha)$. We can then fix some $\epsilon>0$ and a sequence of
times $t_n$ that converges to infinity, such that
$x_{t_n}(\alpha)\leq \tilde y(\alpha)-\epsilon$. In particular,
for any $\delta>0$, we have
\begin{equation}\label{eq:aa}
\int_{\alpha-\delta}^{\alpha} x_{t_n}(\beta)\,d\beta\leq
\int_{\alpha-\delta}^{\alpha} x_{t_n}(\alpha)\,d\beta =\delta
x_{t_n}(\alpha)\leq \delta \tilde y(\alpha)-\delta\epsilon.
\end{equation}

Since $\jnt{\tilde  y} (\cdot)$ is continuous at $\alpha$, we can choose
$\delta$ so that $\jnt{\tilde  y} (\alpha-\delta)\geq \tilde
y(\alpha)-\epsilon/3$. Furthermore, there exists a sequence of
times $\tau_n$ that converges to infinity and such that
$$x_{\tau_n}(\alpha-\delta)\geq \tilde y(\alpha-\delta)-\frac{\epsilon}{3}\geq y(\alpha)-\frac{2\epsilon}{3}.$$
At those times, we have
\begin{equation}\label{eq:bb}
\int_{\alpha-\delta}^{\alpha} x_{\tau_n}(\beta)\, d\beta \geq
\int_{\alpha-\delta}^{\alpha} x_{\tau_n}(\alpha-\delta)\, d\beta
=\delta x_{\tau_n}(\alpha-\delta) \geq \delta \tilde y(\alpha)
-\frac{2\delta\epsilon}{3}.
\end{equation}
However, Eqs.\ \eqref{eq:aa} and \eqref{eq:bb} contradict the
fact that $\int_{\alpha-\delta}^{\alpha} x_t(\alpha)\,d\alpha$
converges, thus establishing the desired result.
\end{proof}

\def\oF{{\overline F}}

We now characterize the fixed points \jnt{and the possible limit
points} of the system. Let $F\subset X$ be the set of
nondecreasing functions $\tilde s$ such that for every
$\alpha,\beta \in I$, \jnt{either $\tilde s(\alpha)  = \tilde
s(\beta)$ or} $\abs{\tilde s(\alpha)-s(\beta)}>1$. \jnt{Similarly,
let $\oF$ be the set of nondecreasing functions $\tilde s$ such
that for almost every pair $(\alpha,\beta)\in I^2$, either $\tilde
s(\alpha)  = \tilde s(\beta)$ or $\abs{\tilde
s(\alpha)-s(\beta)}\geq 1$. } \jnt{Finally, we say} that $\tilde
s\in X$ is a {\em fixed point} if the integral equation
(\ref{eq:def_cont_agent_integral}) with initial condition $\tilde
s$ admits a unique solution $x_t = \tilde s$ for all $t$.

\begin{prop}
$ $
\begin{itemize}
\item[(a)] Let $x$ be a \jnt{nondecreasing} \blue{(in $\alpha$, for all $t$)} solution of the integral equation
(\ref{eq:def_cont_agent_integral}), and suppose that $\tilde
y(\alpha)=\lim_{t\to\infty}x_t(\alpha)$, almost everywhere. Then,
$\tilde y \in \oF$. \item[(b)] If $\tilde s\in F$, then $\tilde
s$ is a fixed point. \item[(c)] If $\tilde s$ is a nondecreasing
fixed point, then $\tilde s\in \oF$.
\end{itemize}
\end{prop}
\begin{proof}
(a) We take the limit in Eq.\ \eqref{eq:inta}, as $t_2\to\infty$.
Since the left-hand side converges, and the integral inside the
brackets is nonnegative (by Eq.\ \eqref{eq:intb}), it follows that
$$\liminf_{\tau\to\infty}
\int_0^c\int_{\jnt{0}}^1 \chi_{x_\tau}(\alpha,\beta)
\prt{x_\tau(\beta)-x_\tau(\alpha)} d\beta\, d\alpha=0.$$ Using
Eq.\ \eqref{eq:intb} and then Fatou's lemma, we obtain
$$
\int_0^c\int_c^1
\liminf_{\tau\to\infty}\chi_{x_\tau}(\alpha,\beta)
\prt{x_\tau(\beta)-x_\tau(\alpha)}\, d\beta\, d\alpha=0.$$ Note
that $x_\tau(\beta)-x_\tau(\alpha)$ converges to $\tilde
y(\beta)-\tilde y(\alpha)$, \jnt{a.e.} If $\chi_{\tilde y}(\alpha,\beta)=1$,
then $\chi_{x_\tau}(\alpha,\beta)=1$ for $\tau$ large enough. This
shows that $\liminf_{\tau\to\infty}\chi_{x_\tau}(\alpha,\beta)\geq
\chi_{\tilde y}(\alpha,\beta)$. We conclude that
$$\int_0^c\int_c^1 \chi_{\tilde y}(\alpha,\beta) (\tilde y(\beta)-\tilde y(\alpha))\,
d\beta\, d\alpha=0.$$ We integrate this equation over all $c\in
[0,1]$, \jnt{interchange the order of integration,} and obtain
$$\int_0^1\int_\alpha^1 \chi_{\tilde y}(\alpha,\beta) (\tilde y(\beta)-\tilde y(\alpha))(\beta-\alpha)\,
d\beta\, d\alpha=0.$$ This implies that for almost \jnt{all pairs}
$(\alpha,\beta)$, \jnt{with $\alpha<\beta$,} (with respect to the two-dimensional Lebesgue
measure), we have $\chi_{\tilde y}(\alpha,\beta) (\tilde
y(\beta)-\tilde y(\alpha))=0$, \jnt{and} $\tilde y(\beta)\geq
\tilde y(\alpha)+1$. This is possible only if $\tilde y\in \oF$
(the details of this last step are elementary and are omitted).

(b) Suppose that $\tilde s\in F$. We have either $\chi_{\tilde
{\jnt{s}}}(\alpha,\beta)=0$, or $\tilde s(\alpha)=\tilde
s(\beta)$. Thus, $\int \chi_{\tilde s}(\alpha,\beta) (\tilde
s(\beta)-\tilde s(\alpha))\, d\beta=0$, for all $\alpha$, and $x_t
= \tilde {\jnt{s}}$ for all $t$ is thus a solution of the system.
We now prove that this solution is unique. \jnt{(Recall that
uniqueness is part of our definition of a fixed point.)}

Since $\tilde s$ is bounded and belongs to $F$, there exists a
positive $\epsilon < 1/2$ such that for all $\alpha,\beta\in I$,
either $\tilde s(\alpha) = \tilde s(\beta)$ or $\abs{\tilde
s(\alpha) - \tilde s(\beta)}> 1 + 3 \epsilon$. Let now $y$ be a
solution of (\ref{eq:def_cont_agent_integral}) with $\tilde s$ as
initial condition.
\jnt{Equation (\ref{eq:def_cont_agent_integral}) readily implies}
that $\abs{y_t(\alpha)-\tilde s(\alpha)} \leq
\epsilon$ for all $t\in [0,\epsilon]$ and $\alpha\in I$.
Therefore, for $t\in [0,\epsilon]$ there holds $\abs{y_t
(\alpha) - y_t(\beta)}<1$ if and only if $\abs{\tilde s(\alpha) -
\tilde s(\beta)}<1$, and $y_t$ is also a solution of the integral
problem
\begin{equation*}
y_t(\alpha) = \tilde s(\alpha) + \int_{\tau = 0}^t
\prt{\int_{\beta:(\alpha,\beta)\in C_{\tilde s}}
\prt{y_\tau(\beta)-y_\tau(\alpha)}d\tau},
\end{equation*}
which unlike (\ref{eq:def_cont_agent_integral}) is a linear system
since $C_{\tilde s}$ is constant. It can be shown, using for
example the Lipschitz continuity of the corresponding linear
operator, that this system admits a
unique solution, so that $y_t = \tilde s$ holds for $t\in
[0,\epsilon]$. Repeating this reasoning, we obtain
$y_t=\tilde s$ for all $t>0$, and $\tilde s$ is thus a fixed
point.

(c) Suppose that $\tilde s$ is a nondecreasing fixed point. By
\jnt{the} definition of \jnt{a} fixed point, the function $x$
defined by $x_t = \tilde s$ for all $t$ is a solution of the
integral equation (\ref{eq:def_cont_agent_integral}). Since it
trivially converges to $\tilde s$ and remains nondecreasing, the
result follows from part (a) of this Proposition.
\end{proof}

The following theorem summarizes the convergence results of this
subsection.

\begin{thm}\label{thm:conv_continuous}
Let $x$ be a solution of the integral equation
(\ref{eq:def_cont_agent_integral}) such that $x_0$ is regular (or,
more generally, such that $x_t$ is nondecreasing for all $t$).
There exists a function $\tilde y \in F$ such that $\lim_{t\to
\infty}
x_t(\alpha) = \jnt{\tilde y}(\alpha)$ holds for almost all $\alpha$.
Moreover, the set of \jnt{nondecreasing} fixed points  contains $F$ and is
contained in $\oF$.
\end{thm}

\subsection{Stability and inter-cluster distances}

As in the discrete case, we call \emph{clusters} the discrete
\jnt{opinion values held by a positive measure set of agents at} a fixed point $\tilde s$. For a cluster $A$, we denote by $W_A$, \jnt{referred to as the}
\emph{weight of the cluster}, the length of the interval $\tilde
s^{-1}(A)$.  By an abuse of language, we also call \jnt{a} cluster the
interval $\tilde s^{-1}(A)$ \jnt{of indices of the associated agents.} In this section, we show that for
regular initial conditions, the limit to which the system
converges satisfies a condition on the inter-cluster distance
similar to \jnt{the one in} Theorem \ref{thm:stab_disc_agent}. From this result, we extract a necessary condition for stability of \jnt{a}
fixed point.

\begin{thm}\label{thm:conv_to_stab}
Let $\tilde x_0\in X$ be an initial opinion function, $x$ the
solution of the integral equation
(\ref{eq:def_cont_agent_integral}), and $\tilde s = \lim_{t\to
\infty} x_t$ the fixed point to which $x$ converges. If $\tilde
x_0$ is regular, then
\begin{equation}\label{eq:inter_clust_dist_cont_nonstrict}
\abs{B-A} \geq 1 +
\frac{\min\{W_A,W_B\}}{\max\{W_A,W_B\}}\end{equation} holds for
any two clusters $A$ and $B$ of $\tilde s$.
\end{thm}
\begin{proof}
The idea of the proof is to rely on the continuity of $x_t$
\jnt{(as a function of $\alpha$)} at each $t$ to guarantee the
presence of perturbing agents between the clusters. \jnt{Then, if
\eqref{eq:inter_clust_dist_cont_nonstrict} is violated, these
perturbing agents will cause a merging of clusters.}

Let $A,B$ be two clusters of $\tilde s$, \jnt{with $A<B$. Since $\tilde s$ is a fixed point, $1\leq B-A$. Let}
$m=\frac{W_AA+W_BB}{W_A+W_B}$ be their center of mass. Condition
(\ref{eq:inter_clust_dist_cont_nonstrict}) is equivalent to
requiring the center of mass to be \jnt{at least unit distance away} from at
least one of the clusters. Suppose, to obtain a contradiction,
that this condition is not satisfied, that is, that $m$ is
\jnt{less than unit distance away}
from \jnt{each of the two} clusters $A$ and $B$.

Since clusters are \jnt{at least one unit apart,} $A$ and $B$ are
necessarily adjacent, and \jnt{since \eqref{eq:inter_clust_dist_cont_nonstrict} is violated, $B-A<2$.} From the monotonicity of $x_t$, \jnt{there exists some}  $c\in I$ such that
$$ \sup \{\alpha: \lim_{t\to \infty}x_t(\alpha) = A\} = c = \inf
\{\alpha: \lim_{t\to \infty}x_t(\alpha) = B\}.$$ \jnt{Moreover,} we have the inclusions
\begin{equation}\label{eq:inclus_convergence}
\begin{array}{lllll}
(c-W_A,c)&\subseteq& \{\alpha: \lim_{t\to \infty}x_t(\alpha) = A\}
&\subseteq&[c-W_A,c]\\
(c,c+W_B)&\subseteq& \{\alpha: \lim_{t\to \infty}x_t(\alpha) = B\}
&\subseteq&[c,c+W_B].
\end{array}
\end{equation}

Let us fix an $\epsilon >0$. Since $x_t(\alpha)$ converges to
$\tilde s(\alpha)$ for almost every $\alpha$, since all $x_t$
are nondecreasing, \jnt{and since clusters are separated by at least one,} there exists a $t'>0$ such that for all $t\geq
t'$, the following implications are satisfied:
\begin{equation}\label{eq:x_close_convergence}
\begin{array}{lllllll}
\alpha &<& c-W_A-\epsilon &\Rightarrow& x_t(\alpha) &\leq& A-1,\\
\alpha &\in& (c-W_A+\epsilon,c-\epsilon) &\Rightarrow& x_t(\alpha) &\in& (A-\epsilon,A+\epsilon),\\
\alpha &\in& (c+\epsilon,c+W_B-\epsilon) &\Rightarrow& x_t(\alpha) &\in& (B-\epsilon,B+\epsilon),\\
\alpha &>& c+W_B+\epsilon &\Rightarrow& x_t(\alpha) &\geq& B+1.
\end{array}
\end{equation}

\providecommand{\lx}[1]{\hat \l_{#1}}

We introduce \jnt{some} new notation. To each function $\tilde x\in X$,
we associate the function $\lx{\tilde x}:\Re \to (-1,1)$ defined
by
\begin{equation*}
\lx{\tilde x}(q) = \int_{\tilde x ^{-1}((q-1,q+1))}\prt{\tilde
x(\beta) - q} d\beta,
\end{equation*}
The value $\lx{\tilde x}(q)$ represents the derivative \jnt{of} the opinion
of an agent \jnt{whose current opinion is} $q$. In particular, the
differential equation (\ref{eq:def_cont_agent_derivative}) can be
rewritten as $\frac{d}{dt}x_t(\alpha) = \lx{x_t}(x_t(\alpha))$.

Let us evaluate  $\lx{x_t}(q)$ for $q\in [B-1 + \epsilon,
A+1-\epsilon]$.
\jnt{(Note that this interval is nonempty, because $B-A<2$.)}
Observe first that $q - 1 \geq A -1 + \epsilon>
A-1$, because $B-A \geq 1$. From the first relation
in (\ref{eq:x_close_convergence}) and the continuity of $x_t$ with
respect to $\alpha$, we obtain
\begin{equation*}
x_t(c-W_A - \epsilon)\leq  A - 1 < q -1.
\end{equation*}
Observe also that $q-1 \leq A \blue{-} \epsilon$. From
the second relation in (\ref{eq:x_close_convergence})
and the continuity
of $x_t$, we obtain
\begin{equation*}
x_t(c-W_A + \epsilon) \geq \blue{A-\epsilon \geq} q -1.
\end{equation*}
A similar argument around $q+1$ shows that
\begin{equation*}
x_t(c+W_b - \epsilon) \leq q +1 < B +1 \leq x_t(c + W_B +
\epsilon).
\end{equation*}
Provided that $\epsilon$ is sufficiently small,  \blue{these inequalities and the monotonicity of $x_t$ imply that}
\begin{equation*}
[c-W_A + \epsilon, c+W_B - \epsilon] \subseteq
x_t^{-1}\prt{(q-1,q+1)} \subseteq [c-W_A -\epsilon, c+W_B +
\epsilon].
\end{equation*}

It also follows from the inclusions (\ref{eq:x_close_convergence})
that $\int_{c-W_A + \epsilon}^{c - \epsilon}
\prt{x_t(\beta) - q} \jnt{d\beta} = W_A ( A - q) +O(\epsilon)$\\
and $\int_{
c + \epsilon}^{c + W_B - \epsilon} \prt{x_t(\beta) - q} \jnt{d\beta} = W_B (
B - q) +O(\epsilon)$. Therefore,
\begin{equation*}
\lx{x_t}(q) = W_A(\jnt{A}-q) + W_B(B-q) + O(\epsilon) = (W_A+W_B) (m-q)
+ O(\epsilon),
\end{equation*}
Observe now that since \jnt{the two} clusters do not satisfy condition
(\ref{eq:inter_clust_dist_cont_nonstrict}), their center of mass
$m$ lies in $(B-1,A+1)$. \jnt{Provided that $\epsilon$ is
sufficiently small, we have} $\jnt{m\in (}B-1+\epsilon,A+1
-\epsilon\jnt{)}$ and therefore, $\lx{x_t}(B-1+\epsilon)
> 0 $ and $\lx{x_t}(A+1-\epsilon) \jnt{<} 0$  for all $t\geq t'$.

\jnt{Recall} that $\tilde x_0\in X_m^M$ for some $m,M>0$. From
Theorem \ref{thm:existence_continuous}, $x$ satisfies the
differential equation (\ref{eq:def_cont_agent_derivative})
$\frac{d}{dt}x_t(\alpha) = \LL(x_t)(\alpha) =
\lx{x_t}(x_t(\alpha))$, and also condition
(\ref{eq:mt<inc_rate<Mt}). \jnt{In particular,} $x_t$ is
continuous and increasing with respect to $\alpha \in I$, for each
$t$. There exists therefore a positive length interval $J$ such
that $x_{t'}(J) \subseteq [B-1+\epsilon,A+1 -\epsilon]$. Since
$\lx{x_t}(B-1+\epsilon)>0$ and $\lx{x_t}(A+1 -\epsilon)<0$ hold
for any $t\geq t'$, and since $\frac{d}{dt}x_t(\alpha) =
\lx{x_t}(x_t(\alpha))$, this implies that $x_t(J) \subseteq
[B-1+\epsilon,A+1 -\epsilon]$ for all $t\geq t'$. \jnt{Since $J$
has positive length, this} contradicts the inclusions
(\ref{eq:inclus_convergence}) on the convergence to the clusters
$A$ and $B$.
\end{proof}

\jnt{We note} that the above proof also applies to any solution
of (\ref{eq:def_cont_agent_integral}) for which $x_t$ is
continuous with respect to $\alpha \in I$, for all $t$.

From \jnt{Theorem \ref{thm:conv_to_stab}}, we can deduce a
necessary condition \jnt{for} the stability of \jnt{a} fixed
point, \jnt{under a classical definition of stability (in contrast
to the nonstandard stability notion introduced for the}
discrete-agent system. Let $\tilde s$ be a fixed point of
(\ref{eq:def_cont_agent_integral}). We say that $\tilde s$ is
stable, if for every $\epsilon>0$ there is a $\delta>0$ such that
if $\norm{\tilde s - \tilde x_0}_1 \leq \delta$, then
$\norm{\tilde s - x_t}_1 \leq \epsilon$ for every $t$ and every
solution $x$ of the integral equation
(\ref{eq:def_cont_agent_integral}) with $\tilde x_0$ as initial
condition. It can be shown that this classical notion of stability
\jnt{is stronger than} the stability with respect to the addition
of a perturbing agent used in Section \ref{sec:discr_stab}. More
precisely, if we view the discrete-agent system as a special case
of the continuum model, stability under the current definition
implies stability with respect to the \jnt{definition} used in
Section \ref{sec:discr_stab}.

\begin{cor}\label{cor:stab_fixed_points_cont}
Let $\tilde s$ be a fixed point of
(\ref{eq:def_cont_agent_integral}). If $\tilde s$ is stable, then
for any two clusters $A$ and $B$,
\begin{equation}\label{eq:cond_stab_fix_points_continuous}
\abs{B-A} \geq 1 + \frac{\min\{W_A,W_B\}}{\max\{W_A,W_B\}}.
\end{equation}
\end{cor}
\begin{proof}
Suppose that $\tilde s$ does not satisfy this condition, and let
$K$ be the infimum of $\norm{\tilde s -\tilde s'}_1$ \jnt{over} all $s'$
satisfying the condition. Clearly, $K>0$. For every $\delta>0$,
there exist $M\geq m
>0$, and $\tilde x_0\in X_m^M$ such that $\norm{\tilde s -  \tilde
x_0}_1 \leq \delta$. Let $x$ be the solution of the integral
equation (\ref{eq:def_cont_agent_integral}) with $\tilde x_0$ as
initial condition, and $s'$ the \jnt{a.e.-limit of} $x_{\jnt{t}}$.
It follows from Theorem
\ref{thm:conv_to_stab} that $s'$ satisfies condition
(\ref{eq:inter_clust_dist_cont_nonstrict}), and therefore that
$\norm{\tilde s - \tilde s'}_1\geq K$.
\blue{Using the dominated convergence theorem, we obtain}
$\lim_{t\to \infty}\norm{\tilde s-x_t}_1 =
\norm{\tilde s-s'}_1$. As a result, $\lim_{t\to
\infty}\norm{\tilde s-x_t}_1 \geq K
>0 $ holds for initial conditions $\tilde x_0$ \jnt{arbitrarily close to
$\tilde s$, and $\tilde s$ is therefore} unstable.
\end{proof}

It is possible to prove that \jnt{the} strict inequality version of
condition (\ref{eq:cond_stab_fix_points_continuous}) is also
necessary for stability.
The proof for the \jnt{case of equality} relies on
modifying the positions of an appropriate
set of agents and
\quotes{creating} some perturbing agents at the weighted average
of the two clusters. See Chapter 10 of \cite{Hendrickx:2008phdthesis} or
Theorem 6 in \cite{BlondelHendrickxTsitsiklis:2009_Krausemodel}
for the same proof applied to Krause's model. We conjecture that
\jnt{the strict inequality version of
condition (\ref{eq:cond_stab_fix_points_continuous})} is also sufficient.

\begin{conj}\label{conj:stab_cont}
A fixed point $\tilde s$ of (\ref{eq:def_cont_agent_integral}) is
stable according to the norm $\norm{\, \cdot \, }_{1}$ if and only
if, for \jnt{any two clusters} $A,B$,
\begin{equation*}
\abs{B-A} > 1 + \frac{\min\{W_A,W_B\}}{\max\{W_A,W_B\}},
\end{equation*}
\end{conj}

\blue{We note that Conjecture \ref{conj:stab_cont} is a fairly
strong statement. It implies, for example,  that multiple clusters
are indeed possible starting from regular initial conditions,
which is an open question at present.}

\section{Relation between \jnt{the discrete and continuum-agent  models}}
\label{sec:link_discrete-continuous}

We now formally establish \jnt{a} connection between the
discrete-agent and the continuum-agent \jnt{models}, and use this
connection to \jnt{argue} that \jnt{the validity of} Conjecture
\ref{conj:stab_cont} \jnt{implies the validity of} Conjecture
\ref{conj:stab_discr_agent}. \jnt{Toward} this purpose, we begin
by proving a result on the continuity of the opinion evolution
with respect to the initial \jnt{conditions.}

\begin{prop}\label{prop:continuity_xt_x0}
Let $x$ be the solution of the  \jnt{continuum model}
(\ref{eq:def_cont_agent_integral}) for some regular initial
condition $\tilde x_0 \in X_m^M$. For every $\epsilon
>0$ and $T>0$, there exists a $\delta>0$ such \blue{if $y$ is a solution of the \jnt{continuum model}
(\ref{eq:def_cont_agent_integral}) and $\norm{ y_0 - \tilde x_0}_{\infty}
\leq \delta$, then} $\norm{y_t-x_t}_{\blue{\infty}} \leq \epsilon$,
for all $t\in [0,T]$.
\end{prop}
\begin{proof}
From Theorem \ref{thm:existence_continuous}, $x_t\in X_{me^{-t}}$
for all $t$.  Lemma \ref{lem:lipscont} then implies that for any
$\tilde y_{\lmod{t}} \in \lmod{Y}$ and any $t\in [0,T]$,
\begin{equation}\label{eq:bound_LLxt-LLy}
\norm{\LL( \blue{\tilde y_t})-\LL(x_t)}_\infty  \leq \prt{2 +
\frac{8}{m}e^t}\norm{ \blue{\tilde y_t}-x_t}_\infty \leq \prt{2 +
\frac{8}{m}e^T}\norm{ \blue{\tilde y_t}-x_t}_\infty.
\end{equation}
For every $\alpha \in I$, we have
\begin{equation*}
y_t(\alpha) -x_t(\alpha) = y_0(\alpha) -\tilde x_0(\alpha) +
\int_{ 0}^t \prt{\LL(y_\tau)(\alpha) -
\LL(x_\tau)(\alpha)}d\tau.
\end{equation*}
It follows from this relation and from the bound
(\ref{eq:bound_LLxt-LLy}) that
\begin{equation*}
\norm{y_t-x_t}_{\infty} - \norm{y_s-x_s}_{\infty} \leq \int_{s}^t \prt{2 + \frac{8}{m}e^T} \norm{y_\tau-x_\tau}_{\infty}d\tau
\end{equation*}
holds for any $0\leq s\leq t\leq T$. This implies that for all
$t\in [0,T]$,
\begin{equation*}
\norm{y_t-x_t}_{\infty} \leq \norm{ y_0-\tilde x_0}_{\infty}
e^{t\prt{2 + \frac{8}{m}e^T}} \leq \norm{ y_0-\tilde x_0}_{\infty}
e^{T\prt{2 + \frac{8}{m}e^T}}.
\end{equation*}
Fix now an $\epsilon >0$ and take $\delta>0$ such that $\delta
e^{T\prt{2 + \frac{8}{m}e^T}} \leq \epsilon$. It follows from the
inequality above that if $\norm{{\jnt{\tilde  y_0}}-\tilde x_0}_{\infty}\leq
\delta $, then $\norm{y_t-x_t}_{\infty} \leq \epsilon $ for every
$t\in [0,T]$.
\end{proof}

The following result shows that \blue{continuum-agent model}
can be interpreted as the limit when $n\to \infty$ of the
\blue{discrete-agent model,} on any time interval of finite length.
To avoid any risk of ambiguity, we use $\xi$ to denote discrete
vectors in the sequel. Moreover, we assume that such vectors are
always sorted (i.e., $j>i\Rightarrow$ $\xi_j\geq \xi_i$). We define
the operator \blue{$G$ that maps a} discrete (nondecreasing)
vector to a function  by $G(\xi)(\alpha) = \xi_{i}$ if
$\alpha \in [\frac{i-1}{n},\frac{i}{n})$, and $G(\xi)(1) =
\xi(n)$, where $n$ is the dimension of the vector $\xi$. Let $\xi$
be a solution of the discrete-agent \jnt{model}
(\ref{eq:def_system_intro}) with initial condition $\xi(0)$. One
can verify that $G(\xi(t))$ is a solution to the \blue{continuum-agent} integral equation
for a  (\ref{eq:def_cont_agent_integral}) with
$G(\xi(0))$ as initial condition. As a result, the \blue{discrete-agent model}  can be simulated by the \blue{continuum-agent model. The next proposition provides  a converse, in some sense, over finite-length time intervals.}

\begin{thm}\label{thm:approx}
Consider a regular initial \blue{opinion} function $\tilde x_0$,
and let $(\xi^{\blue{\brn}})_{n> 0}$ be a sequence of
(nondecreasing) vectors \blue{in} $\Re^n$ such that $\lim_{n\to
\infty}\norm{G(\xi^{\brn}(0)) - \tilde x_0}_\infty = 0$, and such
that for each $n$, \blue{$\xi^{\brn}(0)$ is a proper initial
condition, admitting a unique solution $\xi^{\brn}(t)$. Then,} for
every $T$ and every $\epsilon \blue{>0}$, there exists $n'$ such
that
\begin{equation*}
\norm{G(\xi^{\brn}(t))-x_t}_{\infty}\leq \epsilon
\end{equation*}
holds for all $t\in [0,T]$ and $n\geq n'$.
\end{thm}
\begin{proof}
\blue{The result} follows directly from Proposition
\ref{prop:continuity_xt_x0} and from the fact that
$G(\xi^{\brn}(t))$ is a solution of
(\ref{eq:def_cont_agent_integral}) with the initial condition
$G(\xi^{\brn}(0))$.
\end{proof}

\blue{When $\tilde x$ is regular,} a simple way of building such a sequence
$(\xi^{\brn}(0))_{n> 0}$ is to take $\xi^{\brn}_i \blue{(0)} = \tilde
x_0(i/n)$. Theorem \ref{thm:approx} implies that the
discrete-agent model approximates arbitrarily well the continuum
model for arbitrarily large periods of time, provided that the
initial distribution of discrete opinions approximates
sufficiently well the initial conditions of the continuum
model. \blue{Now recall that according to} Theorem \ref{thm:conv_to_stab}, \blue{and for regular initial conditions, the continuum-agent model}
converges to a fixed point satisfying the
inter-cluster distance condition
(\ref{eq:inter_clust_dist_cont_nonstrict}).  The conjunction of these two results
seems thus to support our Conjecture \ref{conj:stab_discr_agent},
that the discrete-agent model converges to an
equilibrium satisfying this same condition, provided that the
number of \blue{agents} is sufficiently large and that their initial
opinions approximate some regular function. This \blue{argument, however, is incomplete because} the approximation result \blue{in} Theorem
\ref{thm:approx} is only valid \blue{over finite, not infinite, time intervals. Nevertheless, we will} now show that this reasoning would be valid,
with some exceptions, if Conjecture \ref{conj:stab_cont} holds.

\begin{prop}\label{prop:stab_cont=>stab_discr}
\blue{Suppose that $\tilde x_0$ is regular, and suppose that the
limit $\tilde s$ of the resulting solution $x$ of (\ref{eq:def_cont_agent_integral}) is stable and its clusters satisfy
\begin{equation}\label{eq:stab_prop_equiv_conj}
\abs{B-A}>1 + \frac{\min\{W_A,W_B\}}{\max\{W_A,W_B\}}.
\end{equation}
Let $\xi(0)\in \Re^n$ be a
vector whose $n$ entries are randomly and independently selected
according to a probability density function corresponding to
$\tilde x_0$. Then, the clusters of the limit of the corresponding solution
of (\ref{eq:def_system_intro}) satisfy (\ref{eq:def_system_intro}), with probability that tends to 1 as $n\to\infty$.}
\end{prop}

\begin{proof}
\blue{Let $\tilde s=\lim_{t\to\infty} x_t$, which is assumed to be
stable and to satisfy \eqref{eq:stab_prop_equiv_conj}. Since
\eqref{eq:stab_prop_equiv_conj} involves a strict inequality, we
see that there exists some $K>0$ such that the clusters of any
fixed point $s'$ that satisfies $\norm{s'-\tilde s}_1\leq K$ must
also satisfy \eqref{eq:stab_prop_equiv_conj}. Furthermore, since
$\tilde s$ is stable, there exists some $\epsilon>0$ such that if
a solution of the integral equation
(\ref{eq:def_cont_agent_integral}) satisfies $\norm{y_{t'}-\tilde
s }_1 < \epsilon$ for some $t'$, then $\norm{y_{t}-\tilde s }_1
\leq K$ for all $t\geq t'$. To summarize, if a \lmod{converging}
trajectory $y_t$ comes within $\epsilon$ of $\tilde s$, that
trajectory can only converge to a fixed point whose clusters
satisfy \eqref{eq:stab_prop_equiv_conj}.}

\blue{Suppose now that $\xi(0)$ is a vector generated at random,
as in the statement of the proposition, and whose components are
reindexed so that they are nondecreasing. \lmod{It follows from
Kolmogorov-Smirnov theorems (see \cite{Doob:49} for example) and
the regularity of $\tilde x_0$ that} for any given $\delta>0$, the
probability of the event $\norm{G(\xi(0))-\tilde
x_0}_{\lmod{\infty}}<\delta$ converges to 1, as $n\to\infty$.}

\blue{Since $x_t$ converges to $\tilde s$, a.e., the dominated
convergence theorem implies that there exists some $t'$ be such
that $\norm{x_t-\tilde s}_1 < \epsilon/2$. Let now $\xi(t)$ be a
solution of (\ref{eq:def_system_intro}) for the initial condition
$\xi(0)$, the existence of which is guaranteed with probability 1
by Theorem \ref{thm:exist_unique_discr_agents}. Since $G(\xi(t))$
is also solution of the problem (\ref{eq:def_cont_agent_integral})
with initial condition $G(\xi(0))$, Proposition
\ref{prop:continuity_xt_x0} implies that when $\delta$ is chosen
sufficiently small (which happens with high probability when $n$
is sufficiently large), we will have
$\norm{G(\xi(t'))-x_{t'}}_{\infty}<\epsilon/2$, and, consequently,
$\norm{G(\xi(t'))-x_{t'}}_{\lmod{1}}<\epsilon/2$. Therefore, with
probability that tends to 1 as $n$ increases,
\begin{equation*}
\norm{G(\xi(t'))-\tilde s}_1 \leq \norm{G(\xi(t'))-x_{t'}}_1 +
\norm{x_{t'}-\tilde s}_1  < \frac{\epsilon}{2} +
\frac{\epsilon}{2}=\epsilon.
\end{equation*}
It follows that, with probability that tends to 1 as $n$ increases,
the limit $G(\xi(t))$ is a fixed point that satisfies
\eqref{eq:stab_prop_equiv_conj}.}
\end{proof}

\blue{We now use Proposition \ref{prop:stab_cont=>stab_discr} to
establish the connection between our two conjectures. Suppose that
Conjecture \ref{conj:stab_cont} holds. Let $\tilde x_0$ be a
regular initial condition. By Theorem \ref{thm:conv_to_stab}, the
resulting trajectory converges to a fixed point $\tilde s$ that
satisfies the nonstrict inequality
(\ref{eq:inter_clust_dist_cont_nonstrict}). We expect that
generically the inequality will actually be strict, in which case,
according to Conjecture  \ref{conj:stab_cont}, $\tilde s$ is
stable. Therefore, subject to the genericity qualification above,
Proposition \ref{prop:stab_cont=>stab_discr} implies the validity
of Conjecture \ref{conj:stab_discr_agent}.}

\section{Conclusions}

We have analyzed a simple \jnt{continuous-time} multi-agent system
for which the interaction topology depends on the \jnt{agent}
states. \jnt{We worked with the explicit dynamics of the
interaction topology, which raised a number of difficulties, as
the resulting system is highly nonlinear and discontinuous. This
is in contrast to the case of exogenously determined topology
dynamics, which result into time-varying but linear dynamics.}

After \jnt{establishing} convergence to \jnt{a set of} clusters in
which agents share the same opinion, we focused on the
\jnt{inter-cluster} distances. We proposed an explanation
\jnt{for} the experimentally observed distances based on a notion
of stability \jnt{that is tailored to our context.}  This
\jnt{also} led us to conjecture that the probability \jnt{of
convergence} to a stable equilibrium (in which certain minimal
inter-cluster distances are respected), tends to 1 as the number
of agents increases.

We then introduced a variant of \jnt{the model, involving} a
continuum of agents. For regular initial conditions, we proved the
existence and uniqueness of  \jnt{solutions}, the convergence of
the solution to \jnt{a set} of clusters, and a nontrivial bound on
the \jnt{inter-cluster} distances, \jnt{of the same form as} the
necessary and sufficient condition \jnt{stability for the
discrete-agent model.}  Finally, we established \jnt{a} link
between the \jnt{discrete and continuum models}, and proved that
our first conjecture was implied by \jnt{a seemingly} simpler
conjecture.

The results presented here are parallel to, but much stronger than
those that we obtained for Krause's model of opinion dynamics
\cite{BlondelHendrickxTsitsiklis:2009_Krausemodel}. Indeed, we
have provided here a full \jnt{analysis of the continuum model,}
under the mild \jnt{and easily checkable} assumption of regular
initial conditions.

\jnt{The tractability of the model in this paper can be attributed
to (i)  the inherent symmetry of the model, and (ii) the fact that
it runs in continuous time, although} the latter aspect \jnt{also
raised nontrivial} questions related to the existence and
uniqueness of \jnt{solutions.} \jnt{We note however that similar
behaviors have also been} observed for systems \jnt{without such
symmetry.} One can therefore wonder \jnt{whether} the symmetry
\jnt{is really necessary, or just} allows for comparatively
simpler proofs. \jnt{One can similarly wonder whether our results
admit counterparts in models involving high-dimensional opinion
vectors, where one can no longer rely on monotonic opinion
functions and order-preservation results.}

As in our work on Krause's model, our study of the system on a
continuum and the distances between the \jnt{resulting} clusters
uses the fact that the density of agents between \jnt{the clusters
that are being formed is positive at any} finite time. This
however implies that, unlike \jnt{the discrete-agent case,} the
clusters always remain \jnt{indirectly connected}, and it is not
clear \jnt{whether} this permanent connection \jnt{can eventually
force clusters to merge.} \jnt{In} fact, it is \jnt{an open
question whether there exists a regular initial condition that
leads to multiple clusters, although we strongly suspect this to
be the case.} A simple \jnt{proof would consist of} an example of
regular initial conditions \jnt{that admit a closed-form formula
for $x_t$. However, this is difficult because of the discontinuous
dynamics. The only available examples of this type converge to a
single cluster, as for example, in} the case of any two
dimensional distribution of opinions with circular symmetry (see
\cite{CanutoFagnaniTilli:2008}).

\appendix
\section{Existence and uniqueness of solutions to the discrete-agent equation: Proof of Theorem \ref{thm:exist_unique_discr_agents} }\label{appen:exist_sol_discrete}

We sketch here the proof of Theorem
\ref{thm:exist_unique_discr_agents}, \lmod{a full version of which
is available in
\cite{BlondelHendrickxTsitsiklis:2009_proof_ex_unique_discrete}}.
\lmod{Observe first that if $x$ is the unique solution of the
system (\ref{eq:def_system_intro}) for a given initial condition,
then $x_i(t) = x_j(t)$ implies that $x_i(t') = x_j(t')$ holds for
all $t'>t$. Indeed, one could otherwise build another solution by
switching $x_i$ and $x_j$ after the time $t$, in contradiction
with the uniqueness of the solution. Therefore, every initial
condition $\tilde x_0$ satisfying condition (a) automatically
satisfies condition (c).}

Let us now fix the number of agents $n$, and for each graph $G$ on
$n$ vertices, \jnt{with edge set $E$,} define $X_G\subseteq \Re^n$
as the subset in which $\abs{x_i-x_j}<1$ if $(i,j)\in E$, and
$\abs{x_i-x_j}>1$ if $(i,j)\not \in E$. When restricted to $X_G$,
\jnt{the} system (\ref{eq:def_system_intro}) \jnt{becomes} the
linear time invariant \lmod{differential} system
\begin{equation}\label{eq:LTI}
\dot x_i = \sum_{j:(i,j)\in E}(x_j-x_i),
\end{equation}
which admits a unique solution for any initial condition. This
system can be more compactly written as $\dot x = -L_G x$, where
$L_G$ is the Laplacian matrix of the graph $G$.

Consider an initial condition $\tilde x\in \Re^n$, and suppose
that $\jnt{\tilde x \in} X_{G_0}$ for some $G_0$. \jnt{Let} $x_{G_0}$ be the
unique solution of $\dot x= -L_{G_0}x $ \jnt{with} $x_{G_0}(0) = \tilde
x$. If this solution always remains in $X_{G_0}$, it is
necessarily the unique solution of (\ref{eq:def_system_intro}). Otherwise, \jnt{let}
$t_1>0$  be the first time at which $x_{G_0}(t)\in
\partial X_{G_0}$, and set $x(t) = x_{G_0}(t)$ for all $t\in
[0,t_1]$. By the definition of the sets $X_G$, the point $x(t_1)$
also belongs to the \jnt{boundary} of at least one other set
$X_{G_1}$, with $G_1$ and $G_0$ differing only by one edge
$(i,j)$. We consider here the case \jnt{where} $(i,j)\in E_0$,
$(i,j)\not \in E_1$, and $x_i(t_1) > x_j(t_1)$, but a similar
argument can be made in the three other \jnt{possibilities}. We
also \jnt{assume} that $x(t_1)$ belongs to the closure of no other
set $X_G$, and that $(L_{G_0}x(t_1))_i -
(L_{G_\blue{0}}x(t_{\blue{1}}))_j\not = 0$. This assumption does
not always hold, but can be proved to hold \jnt{for all boundary
points \blue{that can be reached,} except for a set that has zero
measure (with the respect to the relative Lebesgue measure defined
on the lower-dinensional boundary).}

Since $x_i(t) -x_j(t) <1$ for $t$ just before $t_1$ and $x_i(t_1)
-x_j(t_1) =1$, there must hold $\lim_{t \jnt{\uparrow}t_1} (\dot x_i(t) - \dot
x_j(t)) \geq 0$, and thus $- (L_{G_0}x(t_1))_i + (L_{G_0}x(t_1))_j
>0$, because we have assumed that \jnt{the latter} quantity \jnt{is nonzero.}
\jnt{Recall} that $G_1$ is obtained from $G_0$ by \jnt{removing the edge}
$(i,j)$. Since $x_i(t_1)-x_j(t_1) =1$, we have
\begin{equation*}
\begin{array}{lll}
-(L_{G_1}x)_i+(L_{G_1}x)_j&=& -(L_{G_0}x)_i+(L_{G_0}x)_j  -
2(x_j(t_1)-x_i(t_1))\\ &=& -(L_{G_0}x)_i+(L_{G_0}x)_j +2\\
&>&0.\end{array}
\end{equation*}
So, if the solution $x$ can be
\jnt{extended} after $t_1$, there must hold $x_i(t) -x_j(t) >1$
for all $t$ in some positive length open interval starting at
$t_1$. This implies that $x\jnt{(t)}\in X_{G_1}$ on some (possibly
smaller) positive length open interval starting at $t_1$, because
$x(t_1)$ is at a positive distance from all sets \jnt{$X_G$ other
than}  $X_{G_0}$ and $X_{G_1}$. On this latter interval, any
solution $x$ must thus satisfy $\dot x = -L_{G_1}x$. This linear
system admits a unique solution $x_{G_1}$ for which $x_{G_1}(t_1)
= x(t_1)$. Moreover, the solution remains in $X_{G_1}$ for some
positive length time period, again because
$-(L_{G_1}x)_i+(L_{G_1}x)_j
> 0$ and because $x(t_1)$ is at a positive distance from all
sets \jnt{$X_G$ other} than $X_{G_0}$ and $X_{G_1}$. If it remains in \jnt{$X_{G_1}$} forever, we \jnt{extend} $x$ by setting $x(t)=x_{G_1}(t)$ on
$[t_1,\infty)$. \jnt{Otherwise, we extend} $x$, \jnt{as before,} on the
interval $[t_1,t_2]$, where $t_2$ is the first time after $t_1$ at
which $x_{G_1}\in
\partial X_{G_1}$. In both cases, $x$ is a solution to
(\ref{eq:def_system_intro}), on $[0,\infty)$ or $[0,t_2]$
respectively, and is \jnt{unique.} Indeed, we have seen that it is
the unique solution on $[0,t_1)$, that any \jnt{extended} solution
should then enter $X_{G_1}$, and that there is a unique solution
entering $X_{G_1}$ at $t_1$ via
$x(t_1)$.

One can prove that, \jnt{for almost all} initial
conditions, this process can be \jnt{continued} recursively without
encountering any  \quotes{problematic boundary points,} \jnt{namely,}  those for which $(L_{G}x)_i -(L_{G}x)_j =0$, or those incident
to more than two sets.\footnote{\jnt{The authors are pleased to acknowledge discussions with Prof.\ Eduardo Sontag on this assertion and its proof.}}  Such a recursive
\jnt{construction ends} after a finite number of transitions if a
solution \jnt{eventually enters and remains forever in a set $X_G$.} In this case, we have
proved the existence of a unique solution
(\ref{eq:def_system_intro}) on $\Re^+$, differentiable everywhere
but on a finite set \jnt{of times}. \jnt{Alternatively, the construction may result in an
infinite sequence of transition times $t_1,t_2,\ldots$.
If this sequence diverges, we have again a unique solution. A problem arises only if this
sequence converges to some finite time $T^*$, in which case, we could only establish existence and uniqueness on
$[0,T^*)$. The following lemma shows that this problematic behavior
will not arise, and concludes the proof of Theorem
\ref{thm:exist_unique_discr_agents}.}

\begin{lem}\label{l:xeno}
Suppose that the \jnt{above recursive construction never encounters problematic boundary points (in the sense defined above), and} produces an infinite sequence of
transition times $t_0,t_1,\ldots$. Then, this sequence diverges,
\jnt{and therefore there exists a unique
solution $x$,  defined for all $t\geq 0$.
Moreover,
$\dot x(t)$ does not converge to 0 when $t\to \infty$.}
\end{lem}
\begin{proof}
Since the sequence $t_1,t_2,\ldots$ of transition times is
infinite,  a nonempty set of agents is involved in an infinite
number of transitions, and there exists a time $T$ after which
every agent involved in a transition will also be involved in a
subsequent one. Consider now a transition occurring at $s_1>T$ and
involving agents $i$ and $j$. We denote by $\dot x_i(s^-_1)$ and
$\dot x_i(s_1^+)$ the limits $\lim_{t \jnt{\uparrow} s_1} \dot
x_i(t)$ and $\lim_{t \jnt{\downarrow s_1}} \dot x_i(t)$
respectively. (Note that these limits exist because away from
boundary points, the function $x$ is continuously differentiable.)

\lmod{Suppose without loss of generality that $x_i>x_j$. If $i$
and $j$ are connected before $s_1$ but not after. \jnt{The update
equation (\ref{eq:def_system_intro}) implies that} $\dot
x_i(s_1^+) = \dot x_i(s_1^-) - (x_j(s_1)-x_i(s_1))$. \jnt{Noting
that} $x_i(s_1) - x_j(s_1) =1 $, \jnt{we conclude that} that $\dot
x_i(s_1^+) = \dot x_i(s_1^-) +1$. Moreover, $x_i-x_j$ must have
been increasing just before $s_1$, so that $\dot x_i(s_1^-) \geq
\dot x_j(s_1^-)$. If on the other hand $i$ and $j$  are connected
after $s_1$ but not before, then $\dot x_j(s_1^+) = \dot
x_j(s_1^-) +1$, and since $x_i-x_j$ must have been decreasing just
before $s_1$, there holds  $\dot x_j(s_1^-) \geq \dot x_i(s_1^-)$.
In either case, there exits an agent $k_1\in \{i,j\}$ for which
$\dot x_{k_1}(s_1^+) = \max\{ \dot x_i(s_1^-),\dot x_j(s_1^-)\}
+1$.} It follows from $s_1>T$ that this agent will get involved in
some other transition at a further time. Call $s_2$ the first such
time.

The definition (\ref{eq:def_system_intro}) of the system implies
that \jnt{in between transitions,} $\abs{\dot x_i(t)}\leq n$ for
all agents. Using (\ref{eq:def_system_intro}) again, this implies
that $\abs{\ddot x_i(t)}\leq 2n ^2$ for all $t$ at which $i$ is
not involved in a transition. Therefore, $\dot x_{k_1}(s_2^-) \geq
\dot x_{k_1}(\jnt{s_1^+}) - 2n^2 (s_2-s_1) = x_i(s_1^-) + 1 - 2n^2
(s_2-s_1)$. Moreover, by the same argument as above, there exists
a $\jnt{k_2}$ for which $\dot x_{\jnt{k_2}} (s_2^+) = \dot
x_{k_1}(s_2^-) + 1 \geq x_i(s_1^-) + 2 - 2n^2 (s_2-s_1)$.
\jnt{Continuing} recursively, we can build an infinite sequence of
transition times $s_1,s_2,\dots$ \jnt{(a} subsequence of
$t_1,t_2,\dots$), such that for every $m$,
$$
\dot x_{k_m}(s_m^+) \geq \dot x_i(s_1^-)  + m - 2n^2(s_m-\jnt{s_1}).
$$
holds for some agent $k_m$. Since all velocities are bounded by
$n$, this implies that $s_m-s_1$ must diverge \jnt{as} $m$ grows,
and therefore that the sequence $t_1,t_2,\dots$ of transition
times \jnt{diverges.}
\end{proof}

Note that the proof of Lemma \ref{l:xeno} provides an explicit
bound on the number of transitions that can take place during any
given time interval.

\section{Existence and Uniqueness solutions to the continuum-agent model: Proof of Theorem \ref{thm:existence_continuous}
}\label{appen:exit_unique_continuous}

\jnt{Let us fix a function $\tilde x_0 \in X_m^M$, with \jnt{$0<m\leq M$.}
Let us also fix some $t_1$ such that
\begin{equation}\label{eq:ineqs}
\frac{m}{2}  \leq m-2Mt \leq me^{-t} \leq M e^{4t/m}
\leq M+\frac{8}{m}t
\leq 2M,
\end{equation}
and
\begin{equation}\label{eq:ineq2}
\prt{2 + \frac{16}{m}} t\leq \frac{1}{2},\qquad e^t\leq 2,
\end{equation}
for all $t\in[0,t_1]$.
We note, for future reference, that $t_1$ can be chosen as a function
$f(m,M)$, where $f$ is continuous and positive.}

Recall that we
defined the operator $G$ \jnt{that maps measurable functions $x:I\times
[0,t_1]\to \Re$ into the set of such functions}  by
\begin{equation*}
(G(x))_t(\alpha) = \jnt{\tilde x}_0(\alpha) + \int_{0}^t
\LL(x_\tau)(\alpha)d\tau.
\end{equation*}
Observe that $x$ is a solution of \jnt{the integral equation}
(\ref{eq:def_cont_agent_integral}) if and only if \jnt{$x_0=\tilde
x_0$ and $x=G(x)$.} Let $P$ be the set of \jnt{measurable}
functions $x:I\times [0,t_1]\to \Re :(\alpha,t)\to x_t(\alpha)$,
such that $x_0 = \tilde x_0$, and such that for all $t\in
[0,t_1]$, we have $x_t \in X_{m-2Mt}^{M+8(M/m)t}$ or in
\jnt{detail,}
\begin{equation}\label{eq:bound_m_M_for_P}
m-2Mt  \leq \frac{ x_t(\beta)- x_t(\alpha)}{\beta-\alpha} \leq
M+8\frac{M}{m}t,
\end{equation}
for all $t\in [0,t_1]$ and all $\beta \not = \alpha$.
\jnt{In particular,} if $x\in P$, then \begin{equation}\label{eq:bound_m_M_loose}
\frac{m}{2} \leq \frac{ x_t(\beta)- x_t(\alpha)}{\beta-\alpha} \leq 2M
\end{equation}
\jnt{for all $t\in [0,t_1]$ and all $\beta \not = \alpha$.}

\jnt{We \blue{note} that $P$, endowed with the $\|\, \cdot \,
\|_{\infty}$ norm, \lmod{defined by} $$\norm{x}_\infty=
\max_{\alpha\in I, t\in [0,t_1]}\abs{x_t(\alpha)},$$ is a
\blue{complete metric space.} We will apply Banach's fixed point
theorem to the operator $G$ on $P$. The first step is to show a
contraction property of $G$.}

\begin{lem}\label{l:contract}
\jnt{The operator} $G$ is contracting on $P$. \jnt{In particular,}
$\norm{G(y)-G(x)}_\infty <
\frac{1}{2}\norm{y-x}_\infty$ for all $x,y\in P$.
\end{lem}

\begin{proof}
\jnt{Let} $x,y\in P$. \jnt{For any $t\in [0,t_1]$, we have $x_t
\in X_{m/2}$ (cf.\ Eq.\ \eqref{eq:bound_m_M_loose}), and} Lemma
\ref{lem:lipscont} implies that
\begin{equation*}
\norm{\LL(y_t)-\LL(x_t)}_\infty \leq
\prt{2+\frac{16}{m}}\norm{y_t-x_t}_\infty.
\end{equation*}
\jnt{Then,} for every $\alpha \in I$,
\begin{eqnarray*}
\abs{(G(y))_t(\alpha) - (G(x))_t(\alpha)} &=&
\abs{\int_0^{t}\prt{\LL(y_\tau)(\alpha)-\LL(x_\tau)(\alpha)} \, d\tau}\\
&\leq& \int_0^{t}\norm{\LL(y_{\jnt{\tau}})-\LL(x_{\jnt{\tau}})}_\infty d\tau\\
&\leq & \int_0^t \prt{2 + \frac{16}{m}} \norm{y_{\jnt{\tau}}-x_{\jnt{\tau}}}_\infty d\tau\\
&\jnt{\leq} & \prt{2 + \frac{16}{m}} t \norm{y-x}_\infty\\  &\leq& \frac{1}{2}
\norm{y-x}_\infty,
\end{eqnarray*}
where the last inequality \jnt{follows from Eq.\ \eqref{eq:ineq2}.}
\end{proof}

\jnt{Before applying Banach's fixed point theorem, we also need to
verify that $G$ maps $P$ into itself.}

\begin{lem}\label{l:into}
If $x\in P$, then $G(x) \in P$.
\end{lem}

\begin{proof}
\jnt{Suppose that} $x\in P$. \jnt{By definition,} $G(x)_0 = \tilde x_0$,
\jnt{and we only} need to prove that $G(x)$ satisfies condition
(\ref{eq:bound_m_M_for_P}). For $t\in [0,t_1]$ and $\alpha\leq \beta$,  we have
\begin{equation*}
G(x)_t(\beta)- G(x)_t(\alpha) = \tilde x_0(\beta) - \tilde
x_0(\alpha) + \int_{0}^t
\prt{\LL(x_\tau)(\beta)-\LL(x_\tau)(\alpha)} d\tau.
\end{equation*}
It follows from \jnt{the first part of} Lemma \ref{lem:muderiv} and from Eq.\
(\ref{eq:bound_m_M_loose}) that
\begin{equation*}
\LL(x_\tau)(\beta)-\LL(x_\tau)(\alpha) \geq -\prt{x_\tau(\beta) -
x_\tau(\alpha)} \geq - 2M (\beta-\alpha).
\end{equation*}
Since $\tilde x_0(\beta) - \tilde x_0(\alpha) \geq
m(\beta-\alpha)$, for any $t\in [0,t_1]$, we have
\begin{equation*}
G(x)_t(\beta)- G(x)_t(\alpha) \geq m(\beta-\alpha) - \int_0^t 2M
(\beta-\alpha)\, d\tau = (m-2Mt)(\beta-\alpha),
\end{equation*}
so that $G(x)$ satisfies the first inequality \jnt{in}
(\ref{eq:bound_m_M_for_P}).

\jnt{We now use the second part of} Lemma
\ref{lem:muderiv} and Eq.\  (\ref{eq:bound_m_M_loose}), \jnt{to obtain}
\begin{equation*}
\LL(x_\tau)(\beta)-\LL(x_\tau)(\alpha) \leq  \frac{4}{m}
\prt{x_\tau(\beta) - x_\tau(\alpha)} \leq 8\frac{M}{m}
(\beta-\alpha).
\end{equation*}
Since $\tilde x_0(\beta) - \tilde x_0(\alpha) \leq
M(\beta-\alpha)$, for any $t\in [0,t_1]$, we have
\begin{equation*}
G(x)_t(\beta)- G(x)_t(\alpha) \leq M(\beta-\alpha) + \int_0^t
8\frac{M}{m} (\beta-\alpha)d\tau =
\Big(M+8\frac{M}{m}t\Big)(\beta-\alpha).
\end{equation*}
Therefore $G(x)$ also satisfies the second inequality in
(\ref{eq:bound_m_M_for_P}), and belongs to $P$.
\end{proof}

\jnt{By Lemmas \ref{l:contract} and \ref{l:into}, $G$ maps $P$
into itself and is a contraction. It follows, from the Banach
fixed point theorem, that there exists some a unique $x^*\in P$
such that $x^*=G(x^*)$. We now show that no other fixed point can
be found outside $P$.}

\begin{lem}\label{l:no_other}
\blue{If a measurable function \lmod{$x:I\times [0,t_1)$}
satisfies} $x=G(x)$, then \blue{it} satisfies condition
(\ref{eq:mt<inc_rate<Mt}) \jnt{and, in particular, $x\in P$.}
\end{lem}

\begin{proof}
\jnt{Suppose that the} function $x:I\times[0,t_1]\to
\Re:(\alpha,t)\to x_t(\alpha)$ \jnt{satisfies} $x=G(x)$, that is,
$x_t(\alpha) = \tilde x_0(\alpha) + \int_0^t\LL(x_\tau)(\alpha)\,
d\tau$ for all $t$ and $\alpha \in I$. It follows from \jnt{the
first part of} Lemma \ref{lem:muderiv} that
\begin{equation*}
\prt{x_t(\beta) - x_t(\alpha)}- \prt{x_{\jnt{0}}(\beta) - x_{\jnt{0}}(\alpha)}
\geq - \int_{\jnt{0}}^t \prt{x_\tau(\beta) -x_\tau (\alpha)}d\tau
\end{equation*}
holds for $\alpha\leq \beta$ and $0\leq s\leq t$. Together with
the fact that $\tilde x_0 \in X_m$, this implies that
\begin{equation*}
x_t(\beta) - x_t(\alpha) \geq e^{-t}\prt{\tilde x_0(\beta) -
\tilde x_0 (\alpha)}\geq me^{-t}(\beta-\alpha),
\end{equation*}
which proves the first inequality \jnt{in} (\ref{eq:mt<inc_rate<Mt}) and
also that $x_t\in X_{me^{-t}}$ for all $t$.
Using this bound, it
follows from \jnt{the second part of} Lemma \ref{lem:muderiv} that
\begin{equation*}\begin{array}{lll}
\prt{x_t(\beta) - x_t(\alpha)} - \prt{x_{\jnt{0}}(\beta) -
x_{\jnt{0}}(\alpha)} &\leq&   \int_{\jnt{0}}^t
\frac{2}{m}e^t\prt{x_\tau(\beta) -x_\tau (\alpha)}d\tau \\&\leq&
\frac{4}{m} \int_{0}^t \prt{x_\tau(\beta) -x_\tau (\alpha)}d\tau,
\end{array}
\end{equation*}
\jnt{where the last inequality follows from \eqref{eq:ineq2}.
Therefore,
\begin{equation*}
x_t(\beta) - x_t(\alpha) \leq
e^{4t/m}(\tilde x_0(\beta)-\tilde x_0(\alpha)) \leq
Me^{4t/m}(\beta-\alpha)
\end{equation*}
where we have used the assumption that $\tilde x_0\in X^M$. This
shows the second inequality in (\ref{eq:bound_m_M_for_P}) and in
particular, that $x\in P$.}
\end{proof}

\jnt{We have shown so far that the integral equation
(\ref{eq:def_cont_agent_integral}) has a unique solution $x^*$
(for $t\in[0,t_1]$), which also belongs to $P$. We argue that it
is also the unique solution to the differential equation
(\ref{eq:def_cont_agent_derivative}). Since} $\LL(x^*_t)$ is
bounded for all $t$, it follows from the equality
\begin{equation}\label{eq:integral_xt}
x_t^* = \tilde x_0 + \int_0^t \LL(x^*_\tau)\, d\tau,
\end{equation}
that $x^*_t$ is continuous with respect to $t$ \jnt{under} the
$\norm{\, \cdot \,}_\infty$ norm, that is, \lmod{there holds}
$\lim_{\tau \to t} \norm{x^*_\tau - x^*_t}_{\infty}=0$ for all
$t$. \jnt{By Lemma \ref{l:no_other}, $x^*$ satisfies} condition
(\ref{eq:mt<inc_rate<Mt}), and then Lemma \ref{lem:lipscont}
implies that $\LL$ is a Lipschitz continuous at every $x^*_t$. The
continuity of $x^*_t$ with respect to $t$ then implies that
$\LL(x^*_t)$ also evolves continuously with $t$. Therefore, we can
differentiate (\ref{eq:integral_xt}), \jnt{to obtain}
$\frac{d}{dt}x^*_t(\alpha) = \LL(x^*_t)(\alpha)$, for all $t$ and
$\alpha$. The function $x^*$ is thus a solution to the
differential \jnt{equation} (\ref{eq:def_cont_agent_derivative}).
Finally, since every solution of the differential \jnt{equation}
is also a solution of the integral \jnt{equation}, which admits a
unique solution, the solution of the differential \jnt{equation}
is also unique.

To complete the proof of the theorem, \jnt{it remains to} show
that the solution $x^*$ can be extended \jnt{to all $t\in\Re$.}
Let $\tilde x_1 = x^*_{t_1}$. It follows from
(\ref{eq:mt<inc_rate<Mt}) that $\tilde x_1 \in X_{m_1}^{M_1}$,
with $m_1 = me^{-t_1}$ and $M_1 = M \jnt{e^{4t/m}}$. \jnt{By
repeating the argument given for $[0,t_1]$, there exists} unique
\jnt{$x^{**}$}, defined on $I\times [t_1,t_2]$, such that
\begin{equation*}
x^{**}_t = \tilde x_1 + \int_{t_1}^t \LL(x^{**}_\tau)\, d\tau,
\end{equation*}
for all $t\in [t_1,t_2]$, \jnt{and where $t_2-t_1=f(m_1,M_1)$.}
One can easily verify that the function obtained by
\jnt{concatenating} $x^*$ and $x^{**}$ is a (unique) solution of
the integral and differential equations on $[0,t_2]$, and that it
satisfies the bound (\ref{eq:mt<inc_rate<Mt}). Repeating this
argument, we show the existence of a unique solution on every
$[0,t_n]$, with
$$t_{n+1} - t_n = \jnt{f(me^{-t_n},M e^{4t_n/m})}.$$
Since this \jnt{recursion} can be written as $t_{n+1} = t_n +
\jnt{g}(t_n)$, with $\jnt{g}$ continuous and positive on
\jnt{$[0,\infty)$, the sequence $t_n$} diverges. \blue{(To see
this, note that if $t_n\leq t^*$ for all $n$, then $g(t_n)\geq
\min_{0 \leq t\leq t^*}g(t)>0$, which proves that $t_n\to\infty$,
a contradiction.)} \jnt{This completes the proof of} Theorem
\ref{thm:existence_continuous}.

\end{document}